\magnification= \magstep1
\input pictex.tex
\input coman.tex
\input amssym.def

\riferimentifuturi

\autobibliografia

\biblitem{Avez}      
        A.~Avez, {\it Limite des quotients pour des marches al\'eatoires 
	sur des groupes}, C.~R.~Acad. Sci.~Paris S\'er.~A {\bf 276} 
        (1973), 317-320. 
 
\biblitem{Avez2}      
        A.~Avez, {\it Entropie des groupes de type fini}, 
        C.~R.~Acad.~Sci.~Paris S\'er.~A {\bf 275} 
        (1972), 1363-1366. 

\biblitem{Bender} 
        E.~A.~Bender, {\it Asymptotic method in enumeration},  
        SIAM review, {\bf 4} (1974), 485-515. 

\biblitem{DB&FZ1} 
        D.~Bertacchi, F.~Zucca, {\it Equidistribution of random walks 
        on spheres}, J.~Stat.~Phys.~{\bf 94} (1999), 91-111. 

\biblitem{DB&FZ2} 
        D.~Bertacchi, F.~Zucca, {\it Uniform asymptotic estimates of 
        transition probabilities on combs}, preprint. 
 
\biblitem{Berta1} D.~Bertacchi, {\it Random walks on Diestel-Leader  
			graphs.}, preprint. 

\biblitem{Brezis}
	H.~Brezis, {\it Analyse Fonctionnelle-Th\'eorie et Applications},
	Masson, Paris, (1983). 

\biblitem{Burioni-Cassi} R.~Burioni, D.~Cassi, 
	{\it Universal properties of spectral dimension}, Phys.~Rev.~Lett. 
	{\bf 76}, 1091-1093 (1996).  

\biblitem{Burioni-Cassi-Vezzani1}  R.~Burioni, D.~Cassi, A.~Vezzani,
	Phys.~Rev.~E {\bf 60}, 1500  (1999); J.~Phys.~A {\bf 32}, 5539 (1999).

\biblitem{Burioni-Cassi-Vezzani} R.~Burioni, D.~Cassi, A.~Vezzani, 
	{\it The type-problem on the average for random walks
	on graphs}, Eur.~Phys.~J.~B {\bf 15}, 665 (2000).

\biblitem{Cartier} P.~Cartier, {\it Fonctions harmoniques sur un arbre}, Symposia 
        Math.~, {\bf 9} (1972) 203-270. 

\biblitem{Cartier2} P.~Cartier, {\it Harmonic analysis on trees},  
	Proc.~Sympos.~Pure~~Math.~, vol.~26, Amer. Math.~Soc., Providence, 
	R.I.,  1972, 419-424. 

\biblitem{Cassi} D.~Cassi, 
	{\it  Local vs Average Behavior on Inhomogeneous Structures: 
	Recurrence on the Average and a Further Extension of Mermin-Wagner 
	Theorem on Graphs}, Phys.~Rev. Lett. {\bf 76}, 2941 (1996).

\biblitem{Cassi-Fabbian} D.~Cassi, L.~Fabbian, {\it
	The spherical model on graphs}, J.~Phys. A {\bf 32},
	L93 (1999). 
 
\biblitem{Cassi-Regina}      
        D.~Cassi, S.~Regina, {\it Random walks on $d$-dimensional 
        comb lattices}, Modern Phys.~Lett. B  {\bf 6} (1992), 1397-1403. 

\biblitem{Cartwright}      
        D.~I.~Cartwright, {\it Some examples of random walks on free 
                products of discrete groups}, Ann.~Mat.~Pura ed Appl.~(IV),  
                {\bf 151} (1988), 1-15. 
  
\biblitem{Cartwright-Soardi} 
        D.~I.~Cartwright, P.~M.~Soardi,  
        {\it Random walks on free products, quotients and amalgams},  
        Nagoya Math.~j. {\bf 102} (1986), 163-180. 
 
\biblitem{Dixmier}      
        J.~Dixmier, {\it Les moyennes invariantes dans les 
        s\'emigroupes et leurs applications} Acta Sci.~Math.~(Sze\-ged),  
        {\bf 12} A (1950), 213-227. 
 
\biblitem{Doob1} 
        J.~L.~Doob, {\it Discrete potential theory and boundaries},
	J.~Math.~Mech.~{\bf 8}, (1959), 433-458. 
 
\biblitem{Dynkin1} 
	E.~B.~Dynkin, {\it The boundary theory of Markov processes 
	(discrete case)}, Uspehi Mat. Nauk {\bf 24} (1969) no. 2 
	(146) 3--42.  
         
\biblitem{Engelking}      
        R.~Engelking, {\it General Topology}, Sigma series in pure
 	mathematics {\bf 6}, Heldermann, Berlin (1989).

\biblitem{ErFePo} P.~Erd\"os, W.~Feller, H.~Pollard, 
        {\it A theorem on power series}, Bull.~Amer.~Soc.~{\bf 55}, 
         (1949) 201-204. 

\biblitem{Fig-Picar} A.~Fig\`a-Talamanca, M.~A.~Picardello, {\it Harmonic
	analysis on free groups}, Lecture Notes in Pure and Appl.~Math.,
	vol.~87, Dekker, New York and Basel, 1987. 

\biblitem{Gerl}      
        P.~Gerl, {\it A local central limit theorem on some groups},  
        in {\it The First Pannonian Symposium on Mathematical Statistics}, 
        Springer Lecture Notes in Statistics, {\bf 8} (1981), 73-82. 
 
\biblitem{Grab-Woess}
        P.~J.~Grabner, W.~Woess, {\it Functional iterations and 
        periodic oscillations for simple random walk on the  
        Sierpi\'nski graph}, Stochastic Processes and Their Applications  
        {\bf 69} (1997), 127-138. 

\biblitem{Guivarc} Y.~Guivarc'h, {\it Sur la loi des grand nombres et le
	rayon spectral d'une marche al\'eatoire}, Ast\'erisque, {\bf
	74} (1980), 47-98. 

\biblitem{Hille}
        E.~Hille, {\it Analytic Function Theory: vol.~I}, Chelsea Publ.~co., 
        New York, N.Y.~(1959).

\biblitem{Ivanov}      
        A.~A.~Ivanov, {\it Bounding the diameter of a distance-regular  
        graph}, Soviet Math.~Dokl. {\bf 28} (1983), 149-152. 

\biblitem{Kanai1} M.~Kanai, {\it Rough isometries and combinatorial 
		approximations of geometries of non-compact Riemannian 
		manifolds.}, J.~Math.~Soc.~Japan {\bf 37} (1985),
		391-413 (4, 6, 9).

\biblitem{Kanai2} M.~Kanai, {\it Rough isometries and the parabolicity
	of Riemannian manifolds.}, J.~Math. Soc.~Japan {\bf 38} 			(1986), 227-238 (1, 4, 9).

\biblitem{Kreyszig} E.~Kreyszig, {\it Introductory functional analysis with
		applications}, John Wiley \& Sons, New York,
		(1978).

\biblitem{Lyon-Sull} T.~Lyons, D.~Sullivan, {\it Function theory,
	random path and covering spaces.},
	J.~Diff. Geom.~{\bf 19}	(1984), 299-323 (4, 7, 9).

\biblitem{Lyons1} T.~Lyons, {\it Instabilty of the Liouville property
	for quasi-isometric Riemannian manifolds and
	reversible Markov Chains.}, J.~Diff.~Geom.~{\bf 26}
	(1987), 33-66 (8, 9).

\biblitem{Lyons} T.~Lyons, {\it A simple criterion for transience of a
        reversible Markov chain.}, Ann.~Prob.~{\bf 11}
        (1983), 393-402.

\biblitem{Macpherson}
        H.~D.~Macpherson, {\it Infinite distance transitive graphs of finite
        valency}, Combinatorica, {\bf 2} (1982), 63-69.

\biblitem{Mark-Mcgu-Thom} St.~Markvorsen, S.~McGuinness, C.~Thomassen,
        {\it Transient random walks on graphs and metric
	spaces, with applications on hyperbolic surfaces.},
	Proc.~London.~Math.~Soc.~{\bf 64} (1992), 1-20 (4, 9).

\biblitem{Nechaev}
         S.~K.~Nechaev, A.~Yu.~Grosberg, A.~M.~Vershik, {\it Random walks on 
         braid groups: Brownian bridges, complexity and statistics}, 
         J.~Phys.~A {\bf 29} (1996),  
         no 10, 2411-2433. 

\biblitem{Pica}      
        M.~Picardello, W.~Woess, {\it Random walks on amalgams}, 
        Monatsh.~Math.~{\bf 100} (1985), 21-33. 

\biblitem{Pic-Woess}      
        M.~Picardello, W.~Woess, {\it The Full Martin Boundary of The 
        Bi-Tree}, Annals of Probability {\bf 22} (1994), No.~4, 2203-2222. 
 
\biblitem{Pic-Woess2}      
        M.~Picardello, W.~Woess, {\it A converse to the mean value property 
        on homogeneous trees}, Trans. of A.M.S.~ {\bf 311} (1999), No.~1,
        209-225. 

\biblitem{Polya} 
        G.~P\'olya, {\it \"Uber eine Aufgabe der 
	Wahrscheinlichkeitstheorie betreffend die Irrfahrt im
	Stra\ss ennetz},  
        Math.~Ann. {\bf 84} (1921), 149-160. 

\biblitem{Rudin1} 
        W.~Rudin, {\it Real and Complex Analysis}, Mc Graw-Hill,  (1987). 
 
\biblitem{Rudin2}
        W.~Rudin, {\it Principles of Mathematical Analysis}, Mc Graw-Hill,  (1953).

\biblitem{Sawyer} 
        S.~Sawyer, {\it Isotropic random walks in a tree},  
        Z.~Wahrsch.~Verw.~Gebiete, {\bf 42} (1978), 279-292. 

\biblitem{Singer1} 
        I.~Singer, {\it Bases in Banach Spaces: I}, Springer-Verlag, Berlin 
        (1970). 
 
\biblitem{Woess3}      
        W.~Woess, {\it Catene di Markov e Teoria del Potenziale nel 
        Discreto}, Quaderno U.M.I. {\bf 41}, Ed.~Pitagora, Bologna (1996). 

\biblitem{Nearest}      
        W.~Woess, {\it Nearest Neighbour Random Walks on Free Product 
        of Discrete Groups}, Bollettino U.~M.~I.~(6) {\bf 5}-B (1986), 
	961-982. 

\biblitem{Survey}      
        W.~Woess, {\it Random walks on infinite graphs and groups - A  
        survey on selected topics}, Bull.~London Math.~Soc. 
	{\bf 26} (1994), 1-60. 

\biblitem{Woess2} 
        W.~Woess, {\it Random walks on infinite graphs and groups},  
        Cambridge Tracts in Mathematics, {\bf 138}, 
	Cambridge Univ. Press, 2000. 

\biblitem{Yamasaki}
        M.~Yamasaki, {\it Discrete potentials on an infinite network},
        Mem.~Fac.~Sci., Shimane Univ. {\bf 13} (1979), 31-44.

\biblitem{Zucca1} F.~Zucca, {\it Mean value properties for harmonic
                        functions on graphs and trees}, to appear on
			Ann.~Mat.~Pura Appl..

\biblitem{tesifabio}
	F.~Zucca, {\it Radiality, asymptotic equidistribution on spheres
		and harmonic functions on infinite graphs}, Ph.~D.~Thesis,
		Univ.~Milano, (1999).


\def\d{\,{\rm d}}

\def\phi{{\varphi}}
\def\epsilon{{\varepsilon}}
\def\eps{{\varepsilon}}
\def\aut{{\rm AUT}}
\def\deg{{\rm deg}}

\def\dom{{\cal D}}
\def\ra{\rightarrow}

\def\Ra{\Rightarrow}

\def\LongRa{\Longleftrightarrow}

\def\sbs{\subseteq}

\def\ttoa{{$\rm TOA_t$}}
\def\troa{{$\rm ROA_t$}}
\def\toa{{\rm TOA}}
\def\roa{{\rm ROA}}

\def\reali{{\Bbb R}}
\def\complessi{{\Bbb C}}
\def\zzz{{\Bbb Z}}
\def\naturali{{\Bbb N}}
\def\boxf{{$\ulcorner \! \! \lrcorner \! \! \! \! \llcorner \! \! \urcorner$}}
\def\QED{{\ifmmode\sq\else{\unskip\nobreak\hfil
\penalty50\hskip1em\null\nobreak\hfil{\boxf}
\parfillskip=0pt\finalhyphendemerits=0\endgraf}\fi} \vskip 12 pt}

\def\parti#1{{{\cal P} \left({#1}\right )}}

\def\infL{{{\rm inf}L}}
\def\supL{{{\rm sup}L}}
\def\Cap{{\rm cap}}

\def\cF{{\cal F}}

\font\fontgrande=cmss17
\centerline{\fontgrande Classification on the average of
random walks}
\vskip 10 pt
{\noindent
\centerline{Daniela Bertacchi ${}^1$
 -- Fabio Zucca ${}^2$}\par\smallskip
\noindent
\centerline{${}^1$ Universit\`a di Milano - Bicocca}\par\noindent
\centerline{Dipartimento di Matematica e Applicazioni}\par\noindent
\centerline{Via Bicocca degli Arcimboldi 8, 20126 Milano, Italy.}\par
\smallskip
\noindent
\centerline{${}^2$ Universit\`a degli Studi di Milano}\par\noindent
\centerline{Dipartimento di Matematica F.~Enriques}\par\noindent
\centerline{Via Saldini 50, 20133 Milano, Italy.}\par}
\vskip 12 pt

{\leftskip=80pt \rightskip=80pt \noindent {\bf Abstract.}\
We introduce a new
method for studying
large scale properties of random walks.
The new concepts of {\it transience} and {\it recurrence
on the average} are compared with the ones introduced in
\cite{Burioni-Cassi-Vezzani} and with the usual ones;
their relationships are analyzed and various examples
are provided. \par
}



\vskip 18pt

\autosez{sec:intro}{Introduction}

Random walks on graphs provide a mathematical model in many scientific
areas, from finance (financial modelling), to physics
(magnetization properties of metals, evolution of gases) and biology
(neural networks, disease spreading).
In particular graphs describe the microscopical structure of solids,
ranging from very regular structures like crystals or ferromagnetic
metals which are viewed as Euclidean lattices, to the irregular
structure of glasses, polymers or biological objects.

Geometrical and physical properties of these discrete structures are
linked by random walks (especially the simple random walk),
which usually describe the diffusion of a particle in these
more or less regular media.

An interesting feature of random walks on graphs is their large
time scale asymptotics which is deeply connected with the concept
of recurrent or transient random walk.
This classification was first introduced by P\'olya for simple random
walks on lattices (see \cite{Polya}) to distinguish between
random walks which return to the starting point with probability one
(these are recurrent), and those whose return probability is less
than one (which are transient).

We observe that in a vertex-transitive graph (such as the lattice $\zzz^d$)
the return probabilities of the simple random walk do not depend on the
starting vertex; but in the case of any irreducible random walk
they may differ from vertex to vertex, although being
strictly less than one
in one vertex is equivalent to being strictly less than one
in any vertex.
The distinction between recurrent and transient random walks is known
as the type-problem (for the type-problem for random walks on infinite
graphs, see \cite{Woess2}).

It has been recently observed that even though the type of a random walk 
describes local properties of the physical model,
average values of return probabilities over all starting sites 
play a key role in the comprehension of the macroscopical behaviour
of the model itself (like spontaneous breaking of continuous
symmetries \cite{Cassi}, critical exponents of the
spherical model \cite{Cassi-Fabbian}, or harmonic
vibrational spectra \cite{Burioni-Cassi}).

These observations lead to the definition of a new type-problem:
the type-problem on the average (see \cite{Burioni-Cassi-Vezzani}).
The definition introduced by Burioni, Cassi and Vezzani is the following:
given the family of the generating functions of the $n$-step return
probabilities of a random walk, $\{F(x,x|z)\}_{x\in X}$, and a reference
vertex $o\in X$ (where $(X,E(X))$ is the graph to which the random walk
is adapted), the random walk is recurrent on the average if 
$$
\lim_{z\ra 1^-}\lim_{n\ra\infty}{\sum_{x\in B(o,n)}F(x,x|z)\over
\vert B(o,n)\vert}=1, \autoeqno{deffis}
$$
and transient on the average if the value of the double limit is less than
$1$ ($B(o,n)$ is the closed ball of center $o$ and radius $n$, $|\cdot|$ denotes
cardinality).

The ``average'' mentioned in the name given to this new type-problem
is a repeated average over balls with fixed center and increasing
radii (of course existence of the limit of these averages is implicitly
required).
This procedure is a particular case of the following: given a
sequence $\{\lambda_n\}_n$ of probability measures
on the set $X$, for each $n$ we consider the average of $F$ with respect
to $\lambda_n$ (that is the expected value of $F$ with respect to 
$\lambda_n$) and then we take the limit of these averages when $n$ goes
to infinity. Note that in definition~\eqref{deffis} one has to evaluate
a further limit (namely the one for $z$ going to $1$) and
$\lambda_{n}(x)=\chi_{B(o,n)}(x)/\vert B(o,n)\vert$.

From a mathematical point of view the definition of this
``limit on the average'' leads to some problems,
like the existence of the limit, the possibility of exchanging the order of
the two limits and the dependence on the reference vertex $o$.

We provide an example of random walk which has no classification on the
average in the above sense
(the simple random walk on a bihomogeneous tree). Thus the
classification on the average is not complete, while the classical one in
recurrent and transient random walks is complete (we will refer to the usual
classification as the ``local'' one, in contrast with the one ``on the
average'').

We then propose a new classification on the average which is complete and is
in many cases an extension of the one given in \cite{Burioni-Cassi-Vezzani}.
For this new definition we analyze its independence on the reference vertex
and provide a sufficient condition which is weaker than the one produced
in \cite{Burioni-Cassi-Vezzani}.
We make comparisons between the former and the new definitions of
classification on the average and with the local one; we study when these
definitions agree and we give examples of random walks which
behave differently according to different classifications
(that is, which are transient with respect to one of these
    classifications and recurrent according to another one).

Another question which naturally arises is what can be said when
averages are taken over general sets (not necessarily balls), that
is when $\{\lambda_n\}_n$ is defined as
$\{\chi_{B_n}/\vert B_n\vert \}_n$ where $\{B_n\}_n$ is an increasing family
of subsets.
Moreover $\{\lambda_n\}_n$ could be a general family of probability
measures (for instance, for some reasons one would like to give
to some subgraphs a greater weight than the one given to other subgraphs).
We deal with these more general averaging procedure
and prove results which generalize the particular cases.

We briefly outline the content of the paper. In Section~\sref{sec:prelim}
we define the {\it limit on the average}, we recall the
distinction between {\it thermodynamical transience and recurrence on the
average} (\ttoa\ and \troa) as defined in \cite{Burioni-Cassi-Vezzani}.

In Section~\sref{average} we show that the simple random walk on
the bihomogeneous tree has no thermodynamical classification
(Examples~\lemmaref{notclass} and \lemmaref{applflux}) and
we introduce our classification on the average of random
walks (\toa\ and \roa).
The rest of the section is devoted to the study of
averages over balls: we prove that under certain conditions
the classification is independent of the centre of the
balls (Proposition~\lemmaref{indip}). We compare the
classical concepts of recurrence and transience with
the corresponding ``on the average'' and the ``thermodynamical''
ones (Theorem~\lemmaref{twodefo} and Examples~\lemmaref{TRT} and \lemmaref{NTD},
see also Table~1).
A flow criterion on the average is stated
(Theorem~\lemmaref{flusso-o}), which should be compared with the
``classical'' one of Lyons~\cite{Lyons} and Yamasaki~\cite{Yamasaki}.

In Section~\sref{sec:subgraphs} we connect the behaviour of the random
walk on the subgraph to the behaviour of the random walk on the whole
graph. Corollary~\lemmaref{subtoa} and Theorem~\lemmaref{subtrano}
deal with the classification on the average, Theorem~\lemmaref{sottografo}
with the thermodynamical one.

Sections~\sref{FBC} and \sref{convexrw} are devoted respectively to
averages over families of finite sets and general averages. The
two appendices present technical results for averages of
general functions and families of power series.

\medskip
 \hskip50pt{\it Table~1: Comparison between the three classifications}
$$\hbox to\hsize{\vbox{\halign{\indent#\hfil\cr
}}\hfill\vbox{\offinterlineskip
\halign{&\vrule#&
  \strut\quad\hfil#\quad\cr
&\omit&\omit&\multispan{8}\hrulefill\cr
\omit&\omit&height2pt&\omit&&\omit&&\omit&&\omit&\cr
\omit&\omit&
&\ \toa,$\,$\ttoa\  && \ \toa,$\,$\troa\ && \ \roa,$\,$\ttoa\ && \ \roa,$\,$\troa\ &\cr
\omit&\omit&height2pt&\omit&&\omit&&\omit&&\omit&\cr
\omit&\omit&height2pt&\omit&&\omit&&\omit&&\omit&\cr
\multispan{11}\hrulefill\cr
height2pt&\omit&&\omit&&\omit&&\omit&&\omit&\cr
height2pt&\omit&&\omit&&\omit&&\omit&&\omit&\cr
& Recurrent && impossible && impossible &&
Ex.~\lemmaref{cubes} && $\zzz^2$ &\cr
&\omit&& (Th.~\lemmaref{twodef}(ii)) && (Th.~\lemmaref{twodef}(ii))
&&\omit&&\omit&\cr
height2pt&\omit&&\omit&&\omit&&\omit&&\omit&\cr
\multispan{11}\hrulefill\cr
height2pt&\omit&&\omit&&\omit&&\omit&&\omit&\cr
height2pt&\omit&&\omit&&\omit&&\omit&&\omit&\cr
& Transient && $\zzz^3$ && impossible && Ex.~\lemmaref{TRT}
&& Ex.~\lemmaref{NTD} &\cr
&\omit&&\omit&& (Th.~\lemmaref{twodef}(iii))
&&\omit&&\omit&\cr
height2pt&\omit&&\omit&&\omit&&\omit&&\omit&\cr
\multispan{11}\hrulefill\cr
}}}
$$

\autosez{sec:prelim}{Basic definitions}

We start giving the general definition of
a large scale average depending on a sequence of probability measures
on an at most countable set $X$ (we will usually think of $X$ as the
vertex set of an infinite, connected and locally finite graph).

\definition{convexlimit}
{Let
$
\lambda=\{\lambda_{n}\}_{n \in \naturali}
$
be a sequence of probability measures on $X$;
we call
{\rm limit on the $\lambda$-average} (or , if there is no ambiguity,
{\rm limit on the average})
the linear map
$$
L_\lambda(f):= \lim_{n \rightarrow +\infty} \sum_{x\in X} f(x)
\lambda_{n}(x).
$$
We call $\dom (L_\lambda)$ the domain of $L_\lambda$, that is
$$
\eqalign{
\dom (L_\lambda):= \Big \{f\in \complessi^X :\, &
\sum_{x\in X} \vert f(x)\vert  \lambda_{n}(x) < +\infty, \, \forall n \in
\naturali,\cr
{} & \hbox{ and } \exists \lim_{n \rightarrow +\infty} \sum_{x\in X} f(x)
\lambda_{n}(x).\Big \}\cr
}
$$
If $A \subseteq X$ is such that $\chi_A \in {\cal D}(L_\lambda)$, then $A$
is called  $L_\lambda$-{\rm measurable} (or briefly {\rm measurable})
and with a slight abuse of notation, we write $L_\lambda(A)$ instead of
$L_\lambda(\chi_A)$ (and we call it the $L_\lambda$-measure of $A$ or
simply the measure of $A$).
}
If ${\cal F}=\{B_n\}_{n\in\naturali}$ is an increasing family of
finite subsets whose union is $X$,
we call {\it limit on the average} with respect to $\cal F$
(we denote it by  $L_{\cal F}$)
the limit on the $\lambda$-average where $\lambda_{n}(x)=
\chi_{B_n}(x)/\vert B_n\vert $.
\smallskip

When $X$ is a metric space (in our case a locally finite, non-oriented
graph with its natural
distance) and $o\in X$,
we 
study the
{\it limit on the average} 
where $\lambda_{n}(x)=\chi_{B(o,n)}(x)/\vert B(o,n)\vert $,
and we will write $L_o$ instead of $L_\lambda$.
\smallskip

The
limits on the average are particular cases of {\it summability methods}
(see for  instance \cite{Kreyszig} Paragraph 4.10); if
$\lim_{n\ra\infty}  \lambda_{n}(x)=0$ for any $x\in X$ (i.e.~every finite
subset of $X$ is  measurable and its measure is zero) then the
limit on the average is called {\it  regular}.
\smallskip

We want now to point out some remarks on the
Definition~\lemmaref{convexlimit}.

\remark{zeroset}
{Given a
limit on the average on $X$, the set of measurable subsets is not,
in general, a
$\sigma$-algebra (nor an algebra: see Proposition~\lemmaref{measurable}).
Anyway it is easy to show that:
(i) if $S$ is a measurable set such that
$L_\lambda(S)=0$ and $S^\prime \subset S$ then $S^\prime$ is measurable too
and $L_\lambda({S^\prime})=0$;
(ii)
if $A$ is measurable and its
measure  is $0$ then for every bounded complex function $f:X \rightarrow
\complessi$, we have that $\chi_A f\in{\cal D}(L_\lambda)$ and
$L_\lambda(\chi_A f)=0$.}

\remark{defreale}
{We defined $L_\lambda$ on complex valued functions mostly for
technical reasons (integrals in the complex field will be needed),
nevertheless we are interested in real valued functions.
In particular functions taking possibly the
values $\pm \infty$ should be admitted (take for instance $f$ equal
to the Green function of a random walk, which we will define in a moment).
To this aim let us consider a function $f:X\rightarrow \reali \cup \{\pm
\infty\}$, such that,
 for all $n \in \naturali$, at least one of the
following conditions holds:
$$
\cases{
\sum_{x\in X:f(x)>0} f(x) \lambda_{n}(x) < +\infty & \cr
\sum_{x\in X:f(x)<0} f(x) \lambda_{n}(x) > -\infty. & \cr}
\autoeqno{dominf}
$$
For every such function
we introduce the {\it upper limit on the $\lambda$-average} and
{\it lower limit on the $\lambda$-average} as
$$\eqalign{
\supL_\lambda (f) &:= \limsup_{n \rightarrow
+\infty} \sum_{x\in X}f(x) \lambda_{n}(x), \cr
\infL_\lambda (f) &:= \liminf_{n \rightarrow
+\infty} \sum_{x\in X} f(x)\lambda_{n}(x) .\cr
}$$
We easily note that if $f$ is any real valued function satisfying the
above condition then
$f\in {\cal D}(L_\lambda)$ if and only if
$\infL_\lambda(f)=\supL_\lambda(f)\in \reali$.}

\remark{infsup1}
{Since any bounded function (hence any characteristic function of a subset of
$X$) satisfies both equations in~\eqref{dominf},
we have that every subset is $\infL_\lambda$-measurable
and $\supL_\lambda$-measurable (note that
these ``measures'' are not even finitely additive, although they
are defined on $\parti{X}$).
Also note that for every $A \subseteq X$, being $A$ $L_\lambda$-measurable with
measure $1$ is equivalent to $\infL_\lambda(A)=1$
(equivalently, $A^c$ is $L_\lambda$-measurable with zero measure).}

Now we define the functions related to random walks, of which we will
consider averages through the main part of this paper.

Given a given random walk $(X,P)$ we denote by
$p^{(n)}(x,y)$ the $n$-step transition
probabilities from $x$ to $y$ ($n\ge0$) and by $f^{(n)}(x,y)$ the probability
that the random walk starting from $x$ hits $y$ for the first time
after $n$ steps ($n\ge 1$). Then we define the
Green function $G(x,y|z)=\sum_{n\ge0}p^{(n)}(x,y) z^n$
and the generating function of the first time return probabilities
$F(x,y|z)=\sum_{n\ge1}f^{(n)}(x,y) z^n$
where $x,y\in X$, $z\in\complessi$
(further details can be found in \cite{Woess2} Chapter 1.B, where $F$ is called $U$).

An irreducible random walk $(X,P)$ is recurrent if $F(x,x):=F(x,x|1)=1$ for some
$x\in X$ (equivalently for all $x$) and transient if $F(x,x)<1$ for some
$x\in X$ (equivalently for all $x$).

We recall here the flow criterion which characterizes transient
networks. One can associate an electric network to
a reversible random walk $(X,P)$ with reversibility measure $m$
in the following way.
We endow any edge with an orientation $e=(e^-,e^+)$ and with a resistance
$r(e)=(m(e^-)p(e^-,e^+))^{-1}$ (in the case of the simple random
walk $r(e)=1$ for every edge $e$).

A flow $u$ from a vertex $x$ to infinity with input $i_0$ is a function
defined on $E(X)$ such that
$$
\sum_{e:e^-=y}u(e)= \sum_{e:e^+=y}u(e)+i_0 \delta_x(y),
\qquad \forall y \in X.
$$
The energy of $u$ is defined as
$<\!u,u\!>:=\sum_{e\in E(X)} u(e)^2r(e)$.
The existence of finite energy flows is related with transience by the
following theorem
(here $\Cap(x)$ is the capacity of the set $\{x\}$: we refer to 
\cite{Woess2} for the definition).

\theorem{flow}
{Let $(X,P)$ be a reversible random walk. The following are equivalent:
\item{(a)} the random walk is (locally) transient;
\item{(b)} there exists $x\in X$ (equivalently for all $x\in X$) such that
it is possible to find a finite energy flow with non-zero input, from
$x$ to infinity;
\item{(c)} there exists $x\in X$ (equivalently for all $x\in X$) such
that $\Cap(x)>0$.
}

We restate the definition of the type-problem according to
\cite{Burioni-Cassi-Vezzani} (and we call it ``thermodynamical'' to
distinguish it from the definition which will be given later).

\definition{avertransient1}
{Let $(X,P)$ be a random walk, and $o\in X$ a
fixed vertex. Suppose that $F(\cdot,\cdot|z) \in {\cal D}(L_o)$, for all
$z \in (\eps,1)$, for some $\eps\in (0,1)$.
The random walk is said {\rm thermodynamically transient on the average}
with respect to $o$ (briefly \ttoa) if
$$
\lim_{z\ra 1^-} L_o(F(z))<1,
\autoeqno{defttoa}
$$
(where $L_o(F(z)):=L_o(F(\cdot,\cdot|z))$),
{\rm thermodynamically recurrent on the average with respect to $o$ (\troa)}
if the limit is equal to $1$.
}

\autosez{average}{The classification on the average (over balls)}

From now on, if not otherwise stated, we will assume that $(X,E(X))$ is
a connected (infinite), locally finite, non oriented
graph, that $o$ is a fixed vertex of $X$, 
that $(X,P)$ is a random walk, not necessarily adapted to the graph $(X,E(X))$,
and that the limit on the average is $L_o$.

Some natural question are:
is any random walk either \ttoa\ or \troa\ (that is, is the
classification on the average complete)? Does the classification depend
on the choice of $o$? Can we reverse the order of the two limits in
equation \eqref{defttoa}?

Regarding the first question, it is not difficult to find examples
of random walks with no thermodynamical classification.

\example{notclass}
{Let us consider the class of bihomogeneous trees (which coincides
with the class of trees which are radial with respect to every point,
see \cite{Zucca1} Proposition~2.9). Despite its property of symmetry,
the simple random walk on a bihomogeneous tree ${\Bbb T}_{n,m}$ (with
$n \not = m$) is neither \troa\ nor \ttoa\
(for the proof, see Example~\lemmaref{applflux}).}

It would be desirable that the classification on the average
would not depend on the choice of the reference vertex $o$.
It has been shown
in \cite{Burioni-Cassi-Vezzani} Section~4 that, if the graph has bounded geometry
and
$$
\lim_{n \rightarrow +\infty} {\vert \partial B(o,n)\vert  \over
\vert B(o,n)\vert } =0
\autoeqno{slowgrowth}
$$
(where $\partial B(o,n):=\{x\in B(o,n): \exists y \not \in B(o,n),
\, (x,y) \in E(X)\}$) for some $o$, then the limit on the average
is independent of the choice of $o$.
This condition is not satisfied, for instance, by any homogeneous tree
of degree greater than $2$, or by any ``fast growing'' graph.

As for the last question, that is whether the limit in equation
\eqref{defttoa} coincides with $L_o(F)$,
in general the answer is no. Anyway,
exploiting the fact that $F$ is a power series with non negative
coefficients one can show that at least when
$\sum_{n\ge 1} k_n$ converges, where $k_n=\sup_{x\in X} f^{(n)}(x,x)$,
then existence of the limit in \eqref{defttoa} implies
existence of $L_o(F)$ and these limits coincide (see
Proposition~\lemmaref{twodefo} (iii) and (iv)).

In order to overcome these difficulties,
we introduce a new classification on the average.

\definition{avertransient}
{Let $(X,P)$ a random walk, and $\{\lambda_{n}\}_n$
a sequence of probabilities measures on $X$, the random walk is called
{\rm transient on the average with respect to
$\lambda$ ($\lambda$-\toa)}  if
$$
\infL_\lambda(F):=\liminf_{n\ra\infty} \sum_{x\in X}
F(x,x)\lambda_n(x)< 1,
$$
{\rm recurrent on the average with respect to
$\lambda$ ($\lambda$-\roa)} if the limit is equal to $1$. }

Since this limit always exists, this classification is complete
(that is, any random walk is either $\lambda$-\toa\
or $\lambda$-\roa).

 In this section, if not otherwise
stated, $\lambda_n=\chi_{B(o,n)}/\vert B(o,n)\vert $ (we
consider classification with $\infL_o$) and we write \toa\
and \roa\ instead of $\lambda$-\toa\ and $\lambda$-\roa.

Now we exhibit a condition implying that this new
classification (which in the rest of this paper we denote
by ``the classification on the average'', in contrast
with the ``thermodynamical'' one defined by \eqref{defttoa})
 does not depend on  the fixed vertex $o$.
As in the thermodynamical case, this condition is a topological
one for the underlying graph.

\proposition{indip}
{Let $(X,E(X))$ be such that there exists $x\in X$ satisfying
$$
\sup_{n\in\naturali} {\vert S(x,n+1)\vert \over\vert B(x,n)\vert }<+\infty,
\autoeqno{eq:indip}
$$
(where $S(x,n+1)$ is the sphere centered in $x$ with radius $n+1$)
then the classification on the average of any random walk
is independent of the choice of $o$.
}
\proof
We note that equation \eqref{eq:indip} holds for
some $x$ if and only if it
holds for any vertex of $X$. It is easy to show that \eqref{eq:indip}
is equivalent, in the case of the $\infL_o$ classification,
to the  requests of Proposition~\lemmaref{twolimits} (i).  \QED

This condition is weaker than the one for the thermodynamical classification
(equation \eqref{slowgrowth}); indeed observe that \eqref{eq:indip} is satisfied by
any graph with bounded  geometry.
On the other hand, bounded geometry is not necessary, as is shown by
the following example.

\example{nongeomlim}
{Given a strictly increasing sequence of natural numbers $\{s_j\}_j$,
such that $s_0\ge1$, and a vertex $x_0$, construct the tree $T$
as follows (see Figure 1, where $s_j=j$).
Each element
on the sphere $S(x_0,m)$ has exactly one neighbour on the sphere
$S(x_0,m+1)$ if $m\neq s_j$ for any $j\in\naturali$ and exactly
$j$ neighbours if $m=s_j$.
If we choose $s_{j+1}\ge s_j+j+1$ then $T$ satisfies equation \eqref{eq:indip}
and has not bounded geometry.}

\vskip 20pt\hskip2truecm
\beginpicture
\setcoordinatesystem units <.4mm,.4mm> point at 0 0
\setlinear
\def\ttt{
\plot 60 35 60 47 /
\multiput {$\scriptscriptstyle\bullet$} at 60 38 *3 0 3 /
\plot 60 47 53.53 71.15 /
\multiput {$\scriptscriptstyle\bullet$} at 60 47 *5 -.776 2.89 /
\plot 60 47 66.47 71.15 /
\multiput {$\scriptscriptstyle\bullet$} at 60 47 *5 .776 2.89 /
\plot 60 47 42.5 64.5 /
\multiput {$\scriptscriptstyle\bullet$} at 60 47 *5 -2.1 2.1 /
\plot 60 47 77.5 64.5 /
\multiput {$\scriptscriptstyle\bullet$} at 60 47 *5 2.1 2.1 /
\plot 56.1 61.48 56.1 71.48 /
\plot 56.1 61.48 51.1 70.08 /
\plot 56.1 61.48 58.9 71.48 /
\plot 56.1 61.48 49.1 68.48 /
\plot 63.9 61.48 63.1 71.48 /
\plot 63.9 61.48 68.9 70.08 /
\plot 63.9 61.48 61.1 71.48 /
\plot 63.9 61.48 70.9 68.48 /
\plot 70.5 57.5 75.5 66.1 /
\plot 70.5 57.5 79.1 62.5 /
\plot 70.5 57.5 72.98 67.1 /
\plot 70.5 57.5 80.1 59.98 /
\plot 49.5 57.5 44.5 66.1 /
\plot 49.5 57.5 40.9 62.5 /
\plot 49.5 57.5 47.02 67.1 /
\plot 49.5 57.5 39.9 59.98 /
}
\def\ttr{
\startrotation by 0 -1 about 60 35
\ttt
\stoprotation
}

\plot 60 35 120 35 /
\multiput{$\scriptscriptstyle\bullet$} at 60 35 *6 10 0 /

\plot 90 5 90 35 /
\multiput{$\scriptscriptstyle\bullet$} at 90 5 *2 0 10 /

\plot 60 35 60 47 /
\multiput {$\scriptscriptstyle\bullet$} at 60 38 *3 0 3 /
\plot 60 47 53.53 71.15 /
\multiput {$\scriptscriptstyle\bullet$} at 60 47 *5 -.776 2.89 /
\plot 60 47 66.47 71.15 /  
\multiput {$\scriptscriptstyle\bullet$} at 60 47 *5 .776 2.89 /
\plot 60 47 42.5 64.5 /  
\multiput {$\scriptscriptstyle\bullet$} at 60 47 *5 -2.1 2.1 /
\plot 60 47 77.5 64.5 /  
\multiput {$\scriptscriptstyle\bullet$} at 60 47 *5 2.1 2.1 /
\plot 56.1 61.48 56.1 71.48 /
\plot 56.1 61.48 51.1 70.08 /
\plot 56.1 61.48 58.9 71.48 /
\plot 56.1 61.48 49.1 68.48 /
\plot 63.9 61.48 63.1 71.48 /
\plot 63.9 61.48 68.9 70.08 /
\plot 63.9 61.48 61.1 71.48 /
\plot 63.9 61.48 70.9 68.48 /
\plot 70.5 57.5 75.5 66.1 /
\plot 70.5 57.5 79.1 62.5 /
\plot 70.5 57.5 72.98 67.1 /
\plot 70.5 57.5 80.1 59.98 /
\plot 49.5 57.5 44.5 66.1 /
\plot 49.5 57.5 40.9 62.5 /
\plot 49.5 57.5 47.02 67.1 /

\startrotation by 0 1 about 60 35 
\plot 60 35 60 47 /
\multiput {$\scriptscriptstyle\bullet$} at 60 38 *3 0 3 /
\plot 60 47 53.53 71.15 /  
\multiput {$\scriptscriptstyle\bullet$} at 60 47 *5 -.776 2.89 /
\plot 60 47 66.47 71.15 /  
\multiput {$\scriptscriptstyle\bullet$} at 60 47 *5 .776 2.89 /
\plot 60 47 42.5 64.5 /  
\multiput {$\scriptscriptstyle\bullet$} at 60 47 *5 -2.1 2.1 /
\plot 60 47 77.5 64.5 /  
\multiput {$\scriptscriptstyle\bullet$} at 60 47 *5 2.1 2.1 /
\plot 56.1 61.48 56.1 71.48 /
\plot 56.1 61.48 51.1 70.08 /
\plot 56.1 61.48 58.9 71.48 /
\plot 56.1 61.48 49.1 68.48 /
\plot 63.9 61.48 63.1 71.48 /
\plot 63.9 61.48 68.9 70.08 /
\plot 63.9 61.48 61.1 71.48 /
\plot 63.9 61.48 70.9 68.48 /
\plot 70.5 57.5 75.5 66.1 /
\plot 70.5 57.5 79.1 62.5 /
\plot 70.5 57.5 72.98 67.1 /
\plot 70.5 57.5 80.1 59.98 /
\plot 49.5 57.5 44.5 66.1 /
\plot 49.5 57.5 40.9 62.5 /
\plot 49.5 57.5 47.02 67.1 /
\stoprotation

\startrotation by -1 0 about 60 35
\plot 60 35 60 47 /
\multiput {$\scriptscriptstyle\bullet$} at 60 38 *3 0 3 /
\plot 60 47 53.53 71.15 /  
\multiput {$\scriptscriptstyle\bullet$} at 60 47 *5 -.776 2.89 /
\plot 60 47 66.47 71.15 /  
\multiput {$\scriptscriptstyle\bullet$} at 60 47 *5 .776 2.89 /
\plot 60 47 42.5 64.5 /  
\multiput {$\scriptscriptstyle\bullet$} at 60 47 *5 -2.1 2.1 /
\plot 60 47 77.5 64.5 /  
\multiput {$\scriptscriptstyle\bullet$} at 60 47 *5 2.1 2.1 /
\plot 56.1 61.48 56.1 71.48 /
\plot 56.1 61.48 51.1 70.08 /
\plot 56.1 61.48 58.9 71.48 /
\plot 56.1 61.48 49.1 68.48 /
\plot 63.9 61.48 63.1 71.48 /
\plot 63.9 61.48 68.9 70.08 /
\plot 63.9 61.48 61.1 71.48 /
\plot 63.9 61.48 70.9 68.48 /
\plot 70.5 57.5 75.5 66.1 /
\plot 70.5 57.5 79.1 62.5 /
\plot 70.5 57.5 72.98 67.1 /
\plot 70.5 57.5 80.1 59.98 /
\plot 49.5 57.5 44.5 66.1 /
\plot 49.5 57.5 40.9 62.5 /
\plot 49.5 57.5 47.02 67.1 /
\stoprotation

\put{\ttt} at 59.2 36

\startrotation by -1 0 about 90 35 
\ttt
\stoprotation

\put{\ttr} at 59.2 27.6

\put{\it Figure~1} [l] at 80 -10

\endpicture
\medskip
Now we start comparing the two classifications on the average and the
local one.

\proposition{twodefo}
{Let $(X,P)$ be a random walk and
let $\infty$ be the point added to $X$ in order to construct
its one point compactification.
\item{(i)} If there exists $A \subseteq X$ measurable, such
that $L_o(A)=1$ and $\lim_{\scriptstyle x \rightarrow \infty\atop
\scriptstyle x\in A} F(x,x) =\alpha$ then
$L_o(F)$ exists, and is equal to $\alpha$. Thus
the random walk is \toa\ (respectively \roa)
if and only if $\alpha < 1$ (respectively $\alpha=1$);
\item{(ii)} if $(X,P)$ is (locally) recurrent then
 $L_o(F)$ exists, is equal to $1$ and the random walk is \roa;
\item{(iii)} if $(X,P)$ is \troa\ then $L_o(F)$ exists,
 is equal to $1$ and the random walk is \roa;
\item{(iv)} if the series $F(x,x)$ is totally convergent (with
respect to $x\in X$) and $(X,P)$ is \ttoa\ then $L_o(F)$ exists,
 it is less than $1$ and the random walk is \toa;
\item{(v)} $(X,P)$ is \roa\ $\Longleftrightarrow$ for every
 $\eps>0$ the set $\{x:F(x,x)\ge 1-\eps\}$ is measurable with measure $1$;
\item{(vi)}  $(X,P)$ is \roa\ $\Longleftrightarrow$ there exists
$A \subseteq X$ measurable, such that $L_o(A)=1$ and
$\lim_{\scriptstyle x \rightarrow \infty\atop
\scriptstyle x\in A} F(x,x)=1$;
\item{(vii)} $(X,P)$ is \toa\ $\Longleftrightarrow$ there exists $A\sbs X$
  such that $\supL_o(A)>0$ and $\sup_A F(x,x)<1$.
\item{}
\vskip -10 pt}

The proof is a particular case of the proof of
Proposition~\lemmaref{twodef}.

Proposition~\lemmaref{twodefo}(iv) states that being
the series $F(x,x)$  totally convergent guarantees that
the classification on the average and the thermodynamical one agree
(if the last one is admissible). Obviously the function
$F(x,x)$ needs not to be totally convergent even in
the case of simple random walks (see Examples~\lemmaref{TRT} and
\lemmaref{cubes}).

Under certain conditions,
the series $F(x,x)$ is indeed totally convergent.

\proposition{ttoaetoa}
{Let $(X,P)$ be a random walk adapted to the graph $(X,E(X))$.
If one of the following conditions holds
then the series $F(x,x)$ is totally convergent (with respect to $x$).
\item{(i)} There exists a subset $\Gamma$ of $\aut(X)$ (the automorphism
group of the graph)
and  a finite subset $X_0 \subset X$
with the property that
for any $y \in X$ there exist $x \in X_0$ and $\gamma \in \Gamma$
such that $\gamma(x)=y$ and $P$ is $\Gamma$-invariant.
\item{(ii)} The radius of convergence of the Green function $G(x,x|z)$
(which is independent of $x$) is $r>1$.
\item{(iii)}  $(X,P)$ is reversible (with reversibility
measure $m$ and total conductance $a(x,y):=m(x)p(x,y)$)
and it satisfies the {\it strong isoperimetric inequality} that is
$$
\sup_{A\subset X
} {m(A) \over s(A)} <+\infty,
$$
where the supremum is taken over finite subsets $A$
and $s(A):=\sum_{x\in A,y\in A^c}
a(x,y)$.\item{}\vskip-19pt}

\proof
We just outline the main points.\hfill\break\noindent
(i) If $y=\gamma(x)$ and $P$ is invariant under the action of $\gamma$ then
$f^{(n)}(x,x)=f^{(n)}(y,y)$. By hypotheses $k_n:=\sup_{x \in X}
f^{(n)}(x,x)=
\max_{x \in X_0} f^{(n)}(x,x) \le\sum_{x \in X_0} f^{(n)}(x,x)$.
Hence   $\sum_{n=0}^\infty k_n \le \sum_{x\in X_0} F(x,x) \le  |X_0|$.
\hfill\break\noindent
(ii) It follows from $f^{(n)}(x,x) \leq p^{(n)}(x,x)
\leq 1/r^n$ which
holds for every $x\in X$ and every $n\in \naturali$.
\hfill\break\noindent
(iii)  See 
\cite{Woess2} Chapter 2 Theorems 10.3 and 10.9
and apply (ii).
\QED

For instance (i) applies to $\Gamma$-invariant random walks, where
$\Gamma$ is a subgroup and $X$ has a finite number of
orbits with respect to $\Gamma$.
This is the case of random walks adapted to Cayley graphs
or of the simple random walk on quasi transitive graphs.

As for condition (ii), an example is given by a locally finite tree
with minimum degree 2 and with
finite upper bound to the lenghts of its unbranched geodesics.

We observe that even if $(X,P)$ is both thermodynamically classifiable
and classifiable on the average, the two classifications may not agree,
as is shown by the following example.

\example{TRT}
{Let $X:=\bigcup_{n\in\naturali} \{n\} \times \zzz_{n+1}$.
For any $n,m\in \naturali$, $p\in \zzz_{n+1},q\in\zzz_{m+1}$,
$(n,p)$ and $(m,q)$ are neighbours
if and only if one of the following holds (see Figure 2)
\item{1)} $p=0_{\zzz_{n+1}}$ and $q=0_{\zzz_{m+1}}$ and $|m-n|=1$,
\item{2)} $m=n$ and $p-q=\pm 1$,
(where $p-q$ is the usual operation in $\zzz_{n+1}$).\hfill\break
\smallskip
\hskip1.5truecm
\beginpicture
\setcoordinatesystem units <.7mm,.7mm> point at 0 0
\setlinear
\plot 60 35 130 35 /
\multiput{$\cdot$} at 132 35 *8 2 0 /
\plot 150 35 170 35 /
\multiput{$\cdot$} at 172 35 *4 2 0 /
\multiput{$\scriptstyle\bullet$} at 60 35 *3 20 0 /
\put{$\scriptstyle\bullet$} at 160 35
\plot 80 35 80 45 /
\put{$\scriptstyle\bullet$} at 80 45
\put{$\scriptscriptstyle (1,1)$} [rt] at 79 44
\plot 100 35 95 43.6 /
\plot 100 35 105 43.6 /
\plot 105 43.6 95 43.6 /
\put{$\scriptstyle\bullet$} at 95 43.6
\put{$\scriptstyle\bullet$} at 105 43.6
\put{$\scriptscriptstyle (2,1)$} [rb] at 97 44.6
\put{$\scriptscriptstyle (2,2)$} [lb] at 103 44.6
\plot 120 35 113 42 /
\plot 120 35 127 42 /
\plot 113 42 120 49 /
\plot 127 42 120 49 /
\put{$\scriptstyle\bullet$} at 127 42
\put{$\scriptstyle\bullet$} at 113 42
\put{$\scriptstyle\bullet$} at 120 49
\put{$\scriptscriptstyle (3,1)$} [rt] at 114 41
\put{$\scriptscriptstyle (3,2)$} [rb] at 122 50
\put{$\scriptscriptstyle (3,3)$} [lt] at 126 41
\plot 160 35 168.6 40 /
\plot 160 35 151.4 40 /
\put{$\scriptstyle\bullet$} at 168.6 40
\put{$\scriptstyle\bullet$} at 151.4 40
\put{$\scriptscriptstyle (n,n)$} [lt] at 168 39
\put{$\scriptscriptstyle (n,1)$} [rt] at 151 39
\setplotarea x from 145 to 170, y from 35 to 60
\setdashpattern <4pt, 4pt>
\setquadratic
\plot 151.4 40 149 44 153 48 /
\plot 168.6 40 171 44 167 48 /
\plot 153 48 160 49 167 48 /
\put{$\scriptscriptstyle (0,0)$} [rt] at 59 34
\put{$\scriptscriptstyle (1,0)$} [rt] at 79 34
\put{$\scriptscriptstyle (2,0)$} [rt] at 99 34
\put{$\scriptscriptstyle (3,0)$} [rt] at 119 34
\put{$\scriptscriptstyle (n,0)$} [rt] at 159 34
\put{\it Figure~2} [l] at 100 25
\endpicture
\medskip
If $\{p_n\}$ is a $(0,1)$-valued sequence such that $p_n^n \uparrow 1$
and $\alpha \in \reali$,
$\alpha <1/3$, then we define the (adapted) transition probabilities
as follows:
$$\matrix{
\hfill p((0,0),(1,0))=&p((1,1),(1,0)):=1, \hfill & \cr
\hfill p((1,0),(1,1)) :=& p_1+(1-p_1)\alpha,\hfill & \cr
p((n,0),(n-1,0)):=&(1-p_n)\alpha,\hfill  & n \geq 1,\hfill \cr
p((n,0),(n+1,0)):=&(1-p_n)(1-2\alpha),\hfill  & n \geq 1,\hfill\cr
p((n,p),(n,p+1)):=&p_n,\hfill  & n \geq 2,\hfill\cr
p((n,p),(n,p-1)):=&(1-p_n),\hfill  & n \geq 2, p\not = 0,\hfill \cr
p((n,0),(n,n-1)):=&(1-p_n)\alpha,\hfill  & n \geq 2.\hfill
}
$$
By using standard stopping time arguments we easily see that this
random walk is locally transient.
\hfill\break
If we denote by $C_n:=\{(n,p):p\in \zzz_{n+1}\}$ for every $n\in \naturali$,
hence for any $x\in C_n$, 
we have that
$f^{(n)}(x,x) \geq p_n^n$ and $f^{(m)}(x,x) \leq 1-f^{(n)}(x,x)$ for all
$m\neq n$.
Thus $\lim_{x \rightarrow \infty} f^{(m)}(x,x)=0$
for any $m\in \naturali$ and if $z \in (0,1)$
by Bounded Convergence Theorem (using $z^m \geq f^{(m)}(x,x) z^m$)
we derive $\lim_{x\rightarrow \infty} F(x,x|z)=0$. Whence for any
regular $\lambda$ we obtain $L_\lambda(F(z)):=L_\lambda(F(\cdot,\cdot|z))=0$
which implies that the random walk is $\lambda$-\ttoa.
\hfill\break
On the other hand $F(x,x) \geq f^{(m)}(x,x)$ for any $x\in X$,
$m\in\naturali$,
hence if $x\in \cup_{m\geq n} C_m$ we have that
$F(x,x) \geq \inf_{m\geq n} p_m^m = p_n^n$ which implies
$\lim_{x \rightarrow \infty} F(x,x) = 1$ and, for any regular
$\lambda$, $L_\lambda(F)=1$
(that is, the random walk is $\lambda$-\roa).
Since the classification on the average and the thermodynamical one
are different this provides an example of a random
walk for which the series $F(x,x)$ is not totally convergent
(Proposition~\lemmaref{twodefo}(iv)).  }

Let us now make some comparisons between the local classification and
the classification on the average of a random walk.
The previous example, which is locally transient, \ttoa\ and \roa, shows
also that while local recurrence imply recurrence
on the average, local transience does not imply transience on the average.
There are also examples of locally transient, \troa\ and
\roa\ random walks, as is shown by the following.

\example{NTD}
{Given the sequence of natural numbers $\{s_j=\sum_{i=1}^j\beta^i\}_{j\ge1}$,
where $\beta\ge2$ is an integer number, $s_0=0$ and $o$ is a vertex,
the construction of the tree $T$ is similar to the one in
Example~\lemmaref{nongeomlim}. Each element
on the sphere $S(o,m)$ has exactly one neighbour on the sphere
$S(o,m+1)$ if $m\neq s_j$ for any $j\ge0$ and exactly
$\alpha$ neighbours if $m=s_j$ ($\alpha\in \naturali$) (Figure 3 represents
the case $\alpha=3$, $\beta=2$).
An application of Theorem~\lemmaref{flow} proves that $T$ is locally transient
if and only if $\alpha>\beta$ (see for instance \cite{Zucca1} 
Remark 4.3).
\hfill\break\indent
It is easy to prove that the set $A$ obtained by removing from $X$
the balls of radius $k$ centered in the elements of $S(o,s_k)$, for all
$k\in \naturali$, has $L_o$-measure equal to $1$.
Moreover on $A$, for every fixed $n$, as $x$ tends to infinity $f^{(n)}(x,x)$
is definitively equal to $f^{(n)}_\zzz(0,0)$ (the first time return
probabilities of the simple random walk on $\zzz$).
Hence by Proposition~\lemmaref{twodefo}(i) the graph is \troa\
(thus \roa)
with respect to any reference vertex.
}
\vskip 20pt\hskip65pt
\beginpicture
\setcoordinatesystem units <.4mm,.4mm> point at 0 0
\setlinear

\def\ttt{
\plot 90 35 90 86 /
\multiput {$\scriptscriptstyle \bullet$} at 90 35 *14 0 3 /

\plot 90 41 61.1 75.5 /
\multiput {$\scriptscriptstyle \bullet$} at 90 41 *12 -1.93 2.29 /
\plot 90 41 118.9 75.5 /
\multiput {$\scriptscriptstyle \bullet$} at 90 41 *12 1.93 2.29 /

\plot 90 53 98.5 84.9 /
\multiput {$\scriptscriptstyle \bullet$} at 90 53 *8 .77 2.9 /
\plot 90 53 81.5 84.9 /
\multiput {$\scriptscriptstyle \bullet$} at 90 53 *8 -.77 2.9 /

\plot 82.3 50.15 68.4 80.05 /
\multiput {$\scriptscriptstyle \bullet$} at 82.3 50.15 *8 -1.27 2.72 /

\plot 97.7 50.15 111.6 80.05 /
\multiput {$\scriptscriptstyle \bullet$} at 97.7 50.15 *8 1.27 2.72 /

\plot 82.3 50.15 55.3 69.05 /
\multiput {$\scriptscriptstyle \bullet$} at 82.3 50.15 *8 -2.46 1.72 /

\plot 97.7 50.15 124.7 69.05 /
\multiput {$\scriptscriptstyle \bullet$} at 97.7 50.15 *8 2.46 1.72 /

\plot 90 77 92.3 85.7 /         
\plot 90 77 87.7 85.7 /
\plot 66.84 68.48 63 76.5 /                     
\plot 113.16 68.48 117 76.5 /

\plot 66.84 68.48 59.5 73.6 /
\plot 113.16 68.48 120.5 73.6 /

\plot 83.84 76.2 83.84 85.2 /
\plot 96.16 76.2 96.16 85.2 /
\plot 83.84 76.2 79.34 83.9 /
\plot 96.16 76.2 100.66 83.9 /

\plot 72.14 71.9 70.06 80.7 /
\plot 107.86 71.9 110 80.7 /
\plot 72.14 71.9 66.4 78.8 /
\plot 107.86 71.9 113.6 78.8 /

\plot 62.62 63.9 56.8 70.8 /
\plot 117.38 63.9 123.2 70.8 /
\plot 62.62 63.9 54.2 67 /
\plot 117.38 63.9 125.8 67 /
}
\ttt

\startrotation by -.5 .86 about 90 35
\ttt
\stoprotation

\startrotation by -.5 -.86 about 90 35
\ttt
\stoprotation

\put{\it Figure~3} [l] at 80 -25
\endpicture
\medskip

It is known that (local) transience is equivalently expressed by one of the
following conditions:
(i) $G(x,x):=G(x,x|1)=+\infty$ for some (i.e.~for every) $x
\in X$;
(ii) $F(x,x)=1$ for some (i.e.~for every) $x \in X$.
In the average case we can only claim a partial result.

\proposition{Ginfinito}
{Let $(X,P)$ be a random walk,.
Then:
\item{(i)} if the random walk is \troa\ then
 $\lim_{z \rightarrow 1^-} \infL_o(G(z))=+\infty$;
\item{(ii)}   if the random walk is \roa\ then
$L_o(G)=+\infty$.
}
For the proof we refer to the general case, see
Proposition~\lemmaref{Ginfinit}.
Observe that in Proposition~\lemmaref{Ginfinito} (i) existence of
$L_o(G(z))$ is not guaranteed and then we have to consider
$\infL_o$ instead. Also notice that reversed implications are not true,
see for instance Example~\lemmaref{hair}
(according to \cite{Burioni-Cassi-Vezzani} this is an example of
a {\it mixed} \ttoa\ graph).

Theorem~\lemmaref{flow} gives a useful tool to (locally) classify reversible random walks.
A similar result can be stated for the classification on the average.

\theorem{flusso-o}
{Let $(X,P)$ be a reversible random walk, with
reversibility
measure $m$ satisfying $\inf m(x)>0$, $\sup m(x)<+\infty$ (in particular
this condition is satisfied by the simple random walk on a graph with
bounded geometry). Then TFAE:
\item{(a)} the random walk is \toa;
\item{(b)} there exists $A\sbs X$ such that $\sup L_o(A)>0$, there
is a finite  energy flow $u^x$ from $x$ to $\infty$ with non-zero input 
$i_0$ for 
every $x\in A$  and $\sup_{x\in A} <\!u^x,u^x\!><+\infty$;  
\item{(c)} there exists $A\sbs X$ such that $\sup L_o(A)>0$ 
 and $\inf_{x\in A} \Cap (x)>0$.}

For the proof see Theorem~\lemmaref{flusso}.

As an application we classify bihomogeneous trees and a whole family
of inhomogeneous trees.
\vskip40pt\hskip20pt
\beginpicture
\setcoordinatesystem units <.5mm,.5mm> point at 0 0
\setlinear
\plot 105 35  105 50  /
\put{$\scriptscriptstyle\bullet$} at 105 35
\put{$\scriptscriptstyle\bullet$} at 105 50
\put{$\scriptscriptstyle\bullet$} at 105 42.5
\plot 105 35  118 27.5  /
\put{$\scriptscriptstyle\bullet$} at 118 27.5
\put{$\scriptscriptstyle\bullet$} at 111.5 31.25
\plot 105 35  92 27.5  /
\put{$\scriptscriptstyle\bullet$} at 92 27.5
\put{$\scriptscriptstyle\bullet$} at 98.5 31.25
\plot 105 50  118 57.5  /
\put{$\scriptscriptstyle\bullet$} at 118 57.5
\put{$\scriptscriptstyle\bullet$} at 111.5 53.75
\plot 105 50  92 57.5  /
\put{$\scriptscriptstyle\bullet$} at 92 57.5
\put{$\scriptscriptstyle\bullet$} at 98.5 53.75
\plot 118 27.5  118 12.5  /
\put{$\scriptscriptstyle\bullet$} at 118 20
\put{$\scriptscriptstyle\bullet$} at 118 12.5
\plot 92 27.5  92 12.5  /
\put{$\scriptscriptstyle\bullet$} at 92 12.5
\put{$\scriptscriptstyle\bullet$} at 92 20
\plot 118 27.5  131 35  /
\put{$\scriptscriptstyle\bullet$} at 131 35
\put{$\scriptscriptstyle\bullet$} at 124.5 31.25
\plot 92 27.5  79 35  /
\put{$\scriptscriptstyle\bullet$} at 79 35
\put{$\scriptscriptstyle\bullet$} at 85.5 31.25
\plot 118 57.5  118 67.5  /
\plot 118 57.5  126.6 52.5  /
\plot 131 35  131 45  /
\plot 131 35  139.6 30  /
\plot 92 57.5  92 67.5  /
\plot 92 57.5  83.4 52.5  /
\plot 79 35  79 45  /
\plot 79 35  70.4 30  /
\plot 92 12.5 83.4 7.5 /
\plot 92 12.5 100.6 7.5 /
\plot 118 12.5 109.4 7.5 /
\plot 118 12.5 126.6 7.5 /
\endpicture

\vskip-80pt\hskip42pt
\beginpicture
\setcoordinatesystem units <.4mm,.4mm> point at 0 0
\setlinear
\put{$\scriptscriptstyle\bullet$} at 90 35

\plot 90 35  90 45  /
\put{$\scriptscriptstyle\bullet$} at 90 45
        \plot 90 45  86 63  /
                \put{$\scriptscriptstyle\bullet$} at 88 54
                \put{$\scriptscriptstyle\bullet$} at 86 63
                        \plot 86 63  82 81  /
                        \put{$\scriptscriptstyle\bullet$} at 84 72
                        \put{$\scriptscriptstyle\bullet$} at 82 81

                        \plot 86 63  70 71  /
                        \put{$\scriptscriptstyle\bullet$} at 70 71
                        \put{$\scriptscriptstyle\bullet$} at 78 67

        \plot 90 45  96 53  /
        \put{$\scriptscriptstyle\bullet$} at 96 53
                \plot 96 53 96 63  /
                \put{$\scriptscriptstyle\bullet$} at 96 63
                        \plot 96 63  91 72 /
                        \put{$\scriptscriptstyle\bullet$} at 91 72

                        \plot 96 63  104 81  /
                        \put{$\scriptscriptstyle\bullet$} at 100 72
                        \put{$\scriptscriptstyle\bullet$} at 104 81

                \plot 96 53  112 65  /
                \put{$\scriptscriptstyle\bullet$} at 104 59
                \put{$\scriptscriptstyle\bullet$} at 112 65

\plot 90 35  108 31  /
\put{$\scriptscriptstyle\bullet$} at 99 33
\put{$\scriptscriptstyle\bullet$} at 108 31
        \plot 108 31  114 39  /
        \put{$\scriptscriptstyle\bullet$} at 114 39
                \plot 114 39  126 55  /
                \put{$\scriptscriptstyle\bullet$} at 120 47
                \put{$\scriptscriptstyle\bullet$} at 126 55

                \plot 114 39  123 35  /
                \put{$\scriptscriptstyle\bullet$} at 123 35

        \plot 108 31 114 23 /
        \put{$\scriptscriptstyle\bullet$} at 114 23
                \plot 114 23  120 31  /
                \put{$\scriptscriptstyle\bullet$} at 120 31
                        \plot 120 31  120 11  /
                        \put{$\scriptscriptstyle\bullet$} at 120 21
                        \put{$\scriptscriptstyle\bullet$} at 120 11

                        \plot 120 31  128 25  /
                        \put{$\scriptscriptstyle\bullet$} at 128 25

                \plot 114 23  98 11  /
                \put{$\scriptscriptstyle\bullet$} at 106 17
                \put{$\scriptscriptstyle\bullet$} at 98 11
                        \plot 98 11  80 19  /
                        \put{$\scriptscriptstyle\bullet$} at 88 15
                        \put{$\scriptscriptstyle\bullet$} at 80 19

                        \plot 98 11  98 1  /
                        \put{$\scriptscriptstyle\bullet$} at 98 1

\plot 90 35  72 31  /
\put{$\scriptscriptstyle\bullet$} at 81 33
\put{$\scriptscriptstyle\bullet$} at 72 31
        \plot 72 31  66 39  /
        \put{$\scriptscriptstyle\bullet$} at 66 39
                \plot 66 39 54 55  /
                \put{$\scriptscriptstyle\bullet$} at 60 47
                \put{$\scriptscriptstyle\bullet$} at 54 55

                \plot 66 39 48 47 /
                \put{$\scriptscriptstyle\bullet$} at 57 43
                \put{$\scriptscriptstyle\bullet$} at 48 47

        \plot 72 31  60 15  /
        \put{$\scriptscriptstyle\bullet$} at 66 23
        \put{$\scriptscriptstyle\bullet$} at 60 15
                \plot 60 15  44 27  /
                \put{$\scriptscriptstyle\bullet$} at 52 21
                \put{$\scriptscriptstyle\bullet$} at 44 27

                \plot 60 15  72 -1  /
                \put{$\scriptscriptstyle\bullet$} at 66 7
                \put{$\scriptscriptstyle\bullet$} at 72 -1

\put {\it Figure~4} [l] at -40 -10
\put {\it Figure~5} [l] at 70 -10
\endpicture
\medskip

\example{applflux}
{Consider the bihomogeneous tree ${\Bbb T}_{m,n}$ and a couple
of vertices $x_n$ and $x_m$, the first with degree $n$ and the second with
degree $m$
(see Figure 4 for the case $m=3$ and $n=2$).
We can construct two finite energy flows $u^n$ and $u^m$
with fixed input $i_0$, respectively from $x_n$ to infinity and from
$x_m$ to infinity.
But then we can obtain a finite energy flux from any vertex $x$ to infinity
(with input $i_0$) by translating $u^n$ or $u^m$ (depending on the degree of
$x$). Thus we can construct a family of fluxes with bounded energy and
this proves that the simple random walk on  ${\Bbb T}_{m,n}$ is \toa. The
proof, which is based on the ideas of Theorem~\lemmaref{subtrano},
can be repeated for any $\lambda$.
Moreover, since $L_o(F)$ does not exist for any reference
vertex $o$, and the series $F(x,x)$ is totally convergent 
by Proposition~\lemmaref{ttoaetoa}(i) (with $X_0=\{x_n,x_m\}$
and $\Gamma=\aut(X)$) then by Proposition~\lemmaref{twodefo}(iii) and (iv)
the simple random walk on  ${\Bbb T}_{m,n}$ cannot be thermodynamically
classifiable.
\hfill\break
Analogously one can show that the simple random walk on a tree $T^\prime_{k,n}$
whose vertices have degree
$2$ or $k$ ($k\ge3$) and such that the distance between ramifications is $n$
($n\ge2$) is \toa\ (while in the former case we had essentially only
two fluxes, here we have at most $[n/2]+1$ fluxes).
\hfill\break
Now consider an inhomogeneous tree $T^{\prime\prime}_{k,n}$
whose vertices have degree $2$ or $k$
($k\ge3$) and such that the distance between ramifications does not exceed
$n$ ($n\ge2$): see Figure 5 for the case $k=3$ and $n=2$. Here the family of finite energy fluxes with fixed input,
constructed from any vertex to infinity, is in general infinite, but the
supremum of the energy is bounded by the supremum of the energy of the
fluxes on $T^\prime_{k,n}$. Hence the simple random walk on
$T^{\prime\prime}_{k,n}$ is \toa.
}

As we have seen, the family of $L_o$-measurable sets plays an important
role in
the classification on the average of random walks, but one has to be careful
when dealing with such sets, since they are not an algebra.
We state the following proposition, as for the proof, see
Proposition~\lemmaref{measurable}.

\proposition{measurableo}
{Let $o\in X$. The class of $L_o$-measurable subsets
 is not an  algebra; in particular there exist two
measurable subsets of $X$, $A$ and $B$, such that $A\cap B$ is not measurable.
  }

\autosez{sec:subgraphs}{Subgraphs and graphs}

In this section we study information which can be inferred from
the knowledge of the behaviour of random walks on subgraphs.

 Since we have to average on a 
subgraph, the first thing to do is to rescale the weights.
We start in a general setting, where $X$ is an at most countable
set, and we average a general function $f$.

\definition{inducedlimit}
{Let $\{\lambda_{n}\}_n$ be a sequence of probability measures
on $X$ and $S \subseteq X$ such that
  $\lambda_{n}(S) >0$, for all $n\in \naturali$. Then
the limit on the average $L_\lambda^S$ defined on $S$ by
$\lambda_{n}^S:=\lambda_{n}|_S / \lambda_{n}(S)$ for
every $n \in \naturali$ and  $x \in S$ is called {\rm
rescaled limit on the average}.
  }

As usual, $L_o^S$ will be the limit in the case of the average over balls.
The following proposition links $L^S_\lambda$  and $L_\lambda$.

\proposition{inducedo}
{Let $S\subseteq X$ be an
$L_\lambda$-measurable subset with positive $L_\lambda$-measure.
If $f\in\complessi^X$, then:
\item{(i)} $f|_S \in {\cal D}(L_\lambda^S)\LongRa
\chi_S \cdot f \in {\cal D}(L_\lambda)$,
\item{(ii)} if $f|_S\in {\cal D}(L_\lambda^S)$ then
${L_\lambda(\chi_S \cdot f) =L_\lambda^S(f|_S) \cdot L_\lambda(S)}.
$
\item{(iii)} ${\infL_\lambda(\chi_S \cdot f)=
\infL_\lambda^S(f|_S)\cdot
L_\lambda(S)}. 
$
}

\proof
The proof is straightforward and we omit it.  \QED

 From now on, we consider the averaging process over balls
of the generating function $F$ related to a random walk $(X,P)$.
There are two different ways of looking at the behaviour of the
random walk on a subgraph $S\sbs X$. The first one is to consider
$S$ as a subset of the graph
with $F(x,x)$ restricted to the
sites in $S$.
The second approach is to view $S$ as an independent graph,
with possibly different generating functions $F(x,x)$.

We start with the first point of view.
Proposition~\lemmaref{inducedo} implies that sets of measure zero have no weight
in the averaging procedure (think of $S$ such that $L_\lambda(S)=1$: then
$L_\lambda^S(F|_S)=L_\lambda(F)$). In
particular, if we can find two subsets such that one of
the two grows strictly slowlier than the other one, then the first one can be
neglected in the averaging process.

\remark{growth1}
{Let $X=X_1\cup X_2$, where $X_1\cap X_2$ is finite. Suppose that
$\vert B(o,n)\cap X_1\vert \le f(n)$, $\vert B(o,n)\cap X_2\vert \ge g(n)$ for
all $n \in \naturali$,
where $f$ and $g$
are two functions such that $f(n)/g(n)$ tends to zero as $n$ goes
to infinity. Then $L_o(X_1)=0$ (hence the classification of any random
walk on $X$ depends only on the restriction of the generating function
$F$ on $X_2$).
}

The following corollary of Proposition~\lemmaref{inducedo} links
classification on a subgraph with classification on the whole graph.

\corollary{subtoa}
{Let $(X,P)$ be a random walk and let $S$ be a subgraph of $X$ such that
$L_o(S)>0$.     If the restriction of $F$ to the subgraph satisfies
 $L_o^S(F|_S)<1$ then $(X,P)$ is \toa.}

  Another result which links the behaviour of
 $F$ on subsets with the behaviour of $F$ on
 the whole graph is the following (note that this result holds for any
 $\lambda$ and any function  $f$ in place of $F$ ).

\proposition{split1}
{Let $\overline X := X \cup \{\infty\}$ be
the one point compactification of $X$ with the discrete topology
and let $\{A_i\}_{i\in \naturali}$ be a partition of
$X$ such that
$A_i$ is $L_o$-measurable and for every $i$ such that
$L_o({A_i})>0$
there exists $\lim_{x \rightarrow \infty} F|_{A_i} (x,x)=:\alpha_i$.
\hfill\break
If  $\sum_{i \in \naturali} F(x,x)
\chi_{A_i}(x)$   is uniformly
convergent with respect to $x\in X$ (to $F$), then
$L_o(F)$ exists and is equal to $\sum_{i=1}^\infty
L_o({A_i}) \alpha_i$ (where $\alpha_i$ can be any real
number if $L_o ({A_i})=0$).
}
\proof
If definitively $A_i=\emptyset$, then the statement follows
by induction on $n$ using
Theorem~\lemmaref{Alexandroff1}
and Proposition~\lemmaref{inducedo}.
\par\noindent
Let us consider the general case.
If $L_o({A_i})>0$ then $\vert A_i\vert =+\infty$ and then $\infty$ is
an accumulation point of $A_i$ in $\overline X$. Since
$A_i \cup \{\infty\}$ with
the induced topology from $\overline X$ is homeomorphic to the one point
compactification of $A_i$ (with the induced topology from $X$), then
it is possible to apply Theorem~\lemmaref{Alexandroff1} to $F|_{A_i}$
obtaining $F|_{A_i} \in {\cal D}(L_o)$. By Proposition~\lemmaref{inducedo}
we have that $\chi_{A_i} F \in {\cal D}(A_i)$ and then
$L_o (\sum_{i=1}^n \chi_{A_i} F) =
\sum_{i=1}^n L_o({A_i}) \alpha_i$.
Using Proposition~\lemmaref{uniformly1} we have the conclusion.
\QED

We remark that even though sets of measure zero have no influence
on the resulting limit on the average of the function $F$ 
their presence may change the return probabilities and hence the function
$F$ 
that we average. This is the main difficulty in the
second approach.

Anyway, under certain regularity conditions we can gain information
on the whole graph from the knowledge of what happens on its
subgraph (regarded as an independent graph).
With the next theorem we give a sufficient condition
for the
simple random walk on a general graph to be \toa\ when  one of
its subgraphs is locally transient (that is the simple random walk
on it is locally transient). The proof is an
application  of Theorem~\lemmaref{flusso-o}.

\theorem{subtrano}
{Let $(A,E(A))$ be a subgraph of $X$ such that
$\supL_o(A)>0$.  Suppose that
there exists $x_0\in A$ such that for every vertex $y \in A$ there
exists an injective map ${\gamma_y} : A \rightarrow A$ such that
(i) ${\gamma_y}(x_0)=y$ and  (ii) for any $w,z \in A$, $(w,z) \in E(A)$
implies
$({\gamma_y}(w), {\gamma_y}(z)) \in E(A)$. If the simple random walk on
$(A,E(A))$
is transient then the simple random walk on $(X,E(X))$ is \toa.
}
The proof will be given in the general case, see Theorem~\lemmaref{subtran}.

We observe that the condition on $A$ in the previous statement is a
requirement of ``self-similarity'' of $A$ (take for instance Cayley graphs).

\corollary{subtran1o}
{Let $(G,E(G))$ be a Cayley graph and
$A \subseteq G$ such that (i) the group identity $e \in A$,
(ii) for any $x,y \in A$ we have
that $xy \in A$ and (iii) the simple random walk on $(A,E(A))$ is transient.
If $(X,E(X))$ is a locally finite graph which contains $(A,E(A))$ as a
subgraph and
$\supL_o(A)>0$ then the simple random walk on $(X,E(X))$ is \toa.
}

We observe that Theorem~\lemmaref{subtrano} and
Corollary~\lemmaref{subtran1o} hold for any $\lambda$ (indeed once
the hypothesis $\supL_\lambda(A)>0$ is satisfied, $\lambda$ plays
no role in the proof).

The last result of this section deals with knowledge of random
walks on subgraphs and the thermodynamical limit on the average.
Before stating it we need a technical lemma and a definition.

\lemma{Xn}
{Let $(X,E(X))$ be a graph with bounded geometry
and $C$ a measurable subset of $X$ such that $L_o(C)=0$.
Let $X_n:=\{x \in X: d(x,C) \leq n\}$,
then $X_n$ is measurable and $L_o(X_n)=0$ for every $n \in \naturali$.}

\proof
We note that for every $n,r \in \naturali$,
$X_n\cap B(o,r) \subseteq \bigcup_{x \in C \cap B(o,n+r)} B(x,n)$ and
by hypotheses,
$$
\vert B(o,r+n)\vert /
\vert B(o,r)\vert \le \vert {B(x,n)}\vert \le M,
$$
where $M=\sup_{x \in X} \vert {B(x,n)} \vert$.
Then
$$
{\vert {X_n\cap B(o,r)}\vert  \over \vert {B(o,r)}\vert } \leq M {\vert
{C\cap B(o,r+n)}\vert
\over \vert {B(o,r)}\vert } \leq M^2 {\vert {C\cap B(o,r+n)}\vert
\over \vert {B(o,r+n)}\vert }\
{\buildrel r \rightarrow +\infty \over \longrightarrow}\  0.
$$
\QED

Our purpose is now to consider a subgraph $(A,E(A))$ as an independent
graph, nevertheless the random walk we study on it should be closely
related to the random walk $(X,P)$  that we suppose adapted
to $(X,E(X))$ (think for instance
of the simple random walk): this is the aim of the following definition.

\definition{indsub}
{Let  $(A,E(A))$ be a subgraph on $(X,E(X))$ and
let $(X,P)$ be a random walk on $(X,E(X))$. A random walk $(A,P_A)$
is called {\rm induced random walk} if for every $x \in A \setminus
\partial A$ and every $y \in A$ we have that
$p(x,y)=p_A(x,y)$.
}
We note that in general the induced random walk is not
uniquely determined, but if $n \in \naturali$ and $x \in A$ are
such that $d(x,\partial A) \geq n$ then
$$
f^{(n)}(x,x)=f_A^{(n)}(x,x), \qquad p^{(n)}(x,x)=p_A^{(n)}(x,x).
\autoeqno{ff}
$$
In the next theorem we deal with a graph which is partitioned in two
subgraphs with known properties. We require that a subset $A$ is
{\it convex} with respect to a vertex $o\in A$, that is that for every $x\in A$
at least one geodesic path from $o$ to $x$ lies in $A$ (hence we are
sure that $d_A(o,x)=d_X(o,x)$ and we denote this distance simply by $d$).

\theorem{sottografo}
{Let $(X,E(X))$ be an infinite graph with bounded
geometry, and let $(A,E(A))$ and $(B,E(B))$ be two subgraphs
such that $\{o\}=A\cap B$, $X=A\cup B$ and $A,B$ are both convex with respect to
$o$. Moreover suppose that
$L_o(A)>0$ and $L_o(\partial A)=0$.
Let $P$ be a stochastic matrix representing a random walk
on $X$ (adapted to $(X,E(X))$) and let us consider two induced random walks
(represented by $P_A$ and $P_B$) on
the subgraphs $(A,E(A))$ and $(B,E(B))$.
Under the previous hypotheses we have that
\item{(i)} any two of the following assertion imply the remaining one:
\itemitem{(i.a)} $(X,P)$ is $L_o$-thermodynamically classifiable;
\itemitem{(i.b)} $(A,P_A)$ is $L_o$-thermodynamically classifiable;
\itemitem{(i.c)} $(B,P_B)$ is $L_o$-thermodynamically classifiable;
\item{(ii)} if two of the assertions in $(i)$ hold then  $(A,P_A)$ \ttoa\
implies $(X,P)$ \ttoa;
\item{(iii)} if two of the assertions in $(i)$ hold
and $L_o(A)<1$ then  $(X,P)$ is \troa\ if and only if $(A,P_A)$ and
$(B,P_B)$ are both \troa.}
\proof
$(i)$ It follows easily by Proposition~\lemmaref{inducedo} and by the
equation~\eqref{sotto1}
below.

$(ii)$ Let $F$ and $F_A$ be the generating functions
(depending on $x\in X$ and $z\in [0,1)$)
of the hitting probabilities
associated to $P$ and $P_A$ respectively.
By the hypotheses $F(z) \in {\cal D}(L_o)$ and
 $F_A(z) \in {\cal D}(L^A_o)$ for all $z\in (\eps,1)$,
 for some $\eps\in (0,1)$.
By equation~\eqref{ff} and Lemma~\lemmaref{Xn} we can
 apply
Proposition~\lemmaref{identity2} to $F_A$ and ${F}_{|A}$ obtaining
 that
${F}_{|A} (z)\in {\cal D}(L_o^A)$ and $L_o^A({F}_{|A}(z))=
L_o^A(F_A(z))$ for all $z\in (\eps,1)$.

From Proposition~\lemmaref{inducedo} we have that $\chi_A F \in {\cal D}
(L_o)$, moreover
$
L_o^A({F_A(z)})\equiv L_o^A({F}_{|A})=L_o({\chi_A F(z)})/ L_o(A).
$
We note that
$$
L_o(F(z))= L_o({\chi_A F(z)})+L_o({\chi_{B}F(z)}) -
L_o({\chi_{\{o\}}F(z)}
) \leq L_o^A({F_A(z)})L_o(A) + (1-L_o(A)),
\autoeqno{sotto1}
$$
hence if $\lim_{z \rightarrow 1^-}L_o^A({F_A(z)}) < 1$ we obviously
have $\lim_{z \rightarrow 1^-} L_o({F(z)})<1$.

$(iii)$ The {\it only if} part is a consequence of (ii).
As for the {\it if} part, let us define $P_*$ by
$$
p_*(x,y)= \cases{p_A(x,y) & if $x,y\in A$, $x \not = o$ \cr
        p_B(x,y) & if $x,y \in B$, $x \not = o$ \cr
        (1/2) p_A(o,y) & if $x=o$, $y \in A$ \cr
        (1/2) p_B(x,o) & if $x=o$, $y \in B$,\cr
        0 & otherwise \cr
}
$$
which is clearly a stochastic matrix satisfying
$p_*(x,y)=p(x,y)$ for every $x \not \in \partial A \cup \partial B$.
Since $\deg(\cdot)$ is bounded, we have that
$L_o(\partial A)=0$
if and only if $L_o(\partial B)=0$;
moreover if we have $x\in X$ and $n \in \naturali$ such that $n \leq d(x,
\partial A \cup  \partial B)$ then
$f_*^{(n)}(x,x)=f^{(n)}(x,x)$. Whence, using again Proposition~
\lemmaref{identity2} (with $C:=\partial A \cup \partial B$
and $X_n$ defined as in Lemma~\lemmaref{Xn}) and
Proposition~\lemmaref{inducedo}, we show that
$L_o({F(z)})  =L_o({F_*(z)})=
L_o^A({F_A(z)})L_o(A)+L_o^B({F_B(z)})L_o(B) \ {\buildrel z
\rightarrow 1^-
 \over \longrightarrow} \ L_o(A)+L_o(B)=1,
$
whence the proof is complete.
\QED

The previous theorem is different from those in \cite{Burioni-Cassi-Vezzani}
since  here a subgraph $A$ is regarded as an independent graph with an
induced random walk. In \cite{Burioni-Cassi-Vezzani},
one is supposed to study the generating
function $F$ of $(X,P)$ to classify the random walk; in our approach
one can study independently two (hopefully) simpler random walk $P_A$ and
$P_B$ (on $A$ and $B$ respectively) and then the classification of the
main random walk can be inferred.

\autosez{FBC}{Averages over increasing sequences of subsets}

In this section we deal with the classification on the average of a random
walk when the average is taken over a family ${\cal F}$ of subsets of $X$.
More precisely, we consider ${\cal F}=\{B_n\}_n$ an increasing sequence
of finite subsets of $X$ such that $\bigcup_nB_n=X$ (we call $\cal F$ an
{\it increasing covering  family} or {\it ICF}),
and denote by $L_{\cal F}$ the corresponding limit on the average.

Clearly, two families of subsets ${\cal F}_1=\{B_n\}_n$
and ${\cal F}_2=\{C_n\}_n$ may give the same classification on the average
of a random walk.
The following proposition provides a sufficient condition for this to
happen.

\proposition{FBC1}
{Given two increasing covering families ${\cal
F}_1=\{B_n\}_n$ and ${\cal F}_2=\{C_n\}_n$ in $X$, such that
\item{(i)} there exist a divergent sequence $\{i_n\}_n$ of natural numbers
and $K>0$ such that $B_{i_n}\supseteq C_n$ and $\vert B_{i_n}\vert /
\vert C_{n}\vert \le K$,
  for every $n$;
\item{(ii)} there exist a divergent sequence $\{j_n\}_n$ of
natural numbers and $K^\prime>0$ such that $C_{j_n}\supseteq B_n$ and
$\vert C_{j_n}\vert /\vert B_{n}\vert \le K^\prime$, for every $n$;
\vskip4pt
\noindent then $L_{{\cal F}_1}$ and $L_{{\cal F}_2}$
induce the same classification on the average of any
random  walk.
}
\proof
These conditions are equivalent to the ones in
Proposition~\lemmaref{twolimits} (i).  \QED

We already observed that the family of $L_o$-measurable sets
is not an algebra.
We now prove it for the more general case of $L_{\cal F}$-measurable sets
(where $\cal F$ is an ICF).

\proposition{measurable}
{Let ${\cal F}=\{B_n\}_n$ be an ICF of $X$.
Then the class of $L_{\cal F}$-measurable subsets
is not an  algebra; in particular there exist $A$, $B$
$L_{\cal F}$-measurable subsets of $X$ such that $A\cap B$
is not  $L_{\cal F}$-measurable.   }
\proof
Let us define, for every $n\in \naturali$,
$m_n:=\vert  B_n\vert $; let $\{A_k, C_k\}_{k \in \naturali}$ be a family of
subsets of  $X$ such that for every $k\in \naturali$, $\{A_k,\, C_k\}$ is a
partition of  $S_k=B_k\setminus B_{k-1}$ with the following two properties
$$
\eqalign{
\vert A_k\vert -\vert C_k\vert &\in \{0, \pm 1\}, \quad \forall k\in \naturali, \cr
\vert \cup_{i=0}^k A_i\vert -\vert \cup_{i=0}^k C_i\vert &  \in \{0, \pm 1\}
, \quad \forall k\in
\naturali,  }
$$
Let us define for every $k \in \naturali$, $a_k:=|\cup_{i=0}^k A_i|$ ,
$c_k:=|\cup_{i=0}^k C_i|$;
since
$X$ is infinite, we can choose an increasing sequence of natural numbers
$\{k_n\}_n$ such that $m_{k_{n+1}}/m_{k_n} \geq 4$.
It is easy to note  that by our hypotheses,
for every $n,i \in \naturali$,
$$
\eqalign{
{1 \over 2} - {1 \over m_n} \leq {a_n \over m_n} &\leq
{1 \over 2} + {1 \over m_n} \cr
{1 \over 2} - {1 \over m_n} \leq {c_n \over m_n} &\leq
{1 \over 2} + {1 \over m_n} \cr
{m_i -2 \over m_n+2} \leq {a_i \over a_n} &\leq
{m_i+2 \over m_n-2}. \cr }
\autoeqno{nuov2}
$$
We finally define the two sets
$$
A:=\bigcup_{i=0}^\infty A_i, \qquad
B:=\bigcup_{i=0}^\infty \left ( \bigcup_{j=k_{2i}+1}^{k_{2i+1}} A_j \cup
\bigcup_{j=k_{2i+1}+1}^{k_{2i+2}} C_j
\right );
$$
by equation~\eqref{nuov2} (since $|A \cap B_n|=a_n$)
we have that $A$ is measurable and $L_{\cF}(A)=1/2$;
similarly
$||B \cap B_n|-c_n | \leq 1+|\{i\in \naturali: k_i <n\}|$.
Since
$m_{k_{n+1}}/m_{k_n} \geq 4$ we have that
$
\lim_{n \rightarrow +\infty} |\{i\in \naturali: k_i <n\}|/ m_n =0
$
(observe that
if $|\{i\in \naturali: k_i <n\}|=j$ then $m_n \geq 4^j$),
then by equation~\eqref{nuov2}
we obtain that $B$ is also measurable and $L_\cF(B)=1/2$.
Moreover
$
A \cap B= \bigcup_{i=0}^\infty \bigcup_{j=k_{2i}+1}^{k_{2i+1}} A_j,
$
hence if $n$ is odd
$$
{|A \cap B \cap B_{k_n}| \over |B_{k_n}|} \geq
{a_{k_n} \over m_{k_n}} \left(1 - {a_{k_{n-1}} \over a_{k_n}}
\right) \  {\buildrel    n \rightarrow +\infty \over \longrightarrow}
\  {1 \over 4};
$$
similarly when $n$ is even (since $m_{k_{n+1}}/m_{k_n} \geq 4$ and
using equation~\eqref{nuov2})
$$
{|A \cap B \cap B_{k_n}| \over |B_{k_n}|}
\leq {a_{k_{n-1}} \over m_{k_n}}.
$$
Since $\displaystyle
\liminf_{n \rightarrow +\infty} {a_{k_{n-1}} \over m_{k_{n-1}}}
{m_{k_{n-1}} \over m_{k_n}} \leq 1/8$,
we have that
$\infL_\cF(A \cap B) \leq 1/8$ meanwhile $1/4 \leq
\supL_\cF(A \cap B)$
which implies that $A \cap B$ is not measurable.
\QED

We give two examples of averages over subsets which are not balls.
The first one shows that with two different ICFs the classification of
a random walk can be different; it is also an example of a random
walk which is \toa\ even if $L_o(G)=+\infty$ (recall
Proposition~\lemmaref{Ginfinito}).
The second one is an example of classification on the average with
an ICF which appears natural and with respect to which the random
walk is \roa\ and \ttoa.

\example{hair}
{Let $X$ be the graph obtained from $\zzz^3$ by deleting all horizontal
edges joining vertices with positive height (compare with
\cite{Burioni-Cassi-Vezzani}
where this graph is an example of mixed \ttoa\ and see Figure 6): we call $X_+$ the set of
vertices with (strictly) positive height and $X_-=X_+^c$.  The simple random
walk on $X$ is
locally transient; using Theorem~\lemmaref{flow}
one can construct a finite energy flow $u$
defined on $E(\zzz^3)$ from the origin $o$ to $\infty$ with input $1$.
By Corollary~\lemmaref{subtran1o} we have that the simple random walk
is \toa\ (more generally it is \toa\ with respect to any $\lambda$
such that $\supL_\lambda(X_-)>0$).
Note that $L_o(G)=\infty$ (more generally given a general limit
in the average $L_\lambda(G)=+\infty$ if and only if $L_\lambda(X_-)>0$).
\hfill\break\indent
The same graph is \troa\ (thus \roa)
with respect to the following ICF: ${\cal F}=\{B_n=(B(o,2^n)\cap X_+)
\cup (B(o,n)\cap X_-)\}_n$.
In this case $X_-$ has $L_\cF$-measure zero and by Lemma~\lemmaref{Xn}, for
every $n$, $f^{(n)}(x,x)=f^{(n)}_\zzz$ outside a set of $L_\cF$-measure zero.
Thus by Theorem~\lemmaref{Alexandroff1} every $f^{(n)}$ has $L_\cF$ limit
equal to $f^{(n)}_\zzz$ and the graph is \troa\
(Theorem~\lemmaref{limit1}(a.1)). }

$\ $
\vskip 60pt\hskip20pt
\beginpicture
\setcoordinatesystem units <.3mm,.3mm> point at 0 0
\setlinear
\plot 90 8 90 72 /
\plot 80 8 80 72 /
\plot 70 8 70 72 /
\plot 100 8 100 72 /
\plot 110 8 110 72 /
\plot 120 8 120 72 /
\plot 130 8 130 72 /

\plot 63 15 63 79 /
\plot 56 22 56 86 /
\plot 49 29 49 93 /

\plot 73 57 73 79 /
\plot 83 57 83 79 /
\plot 93 57 93 79 /
\plot 103 57 103 79 /
\plot 113 57 113 79 /
\plot 123 57 123 79 /
\plot 66 64 66 86 /
\plot 76 64 76 86 /
\plot 86 64 86 86 /
\plot 96 64 96 86 /
\plot 106 64 106 86 /
\plot 116 64 116 86 /
\plot 59 71 59 93 /
\plot 69 71 69 93 /
\plot 79 71 79 93 /
\plot 89 71 89 93 /
\plot 99 71 99 93 /
\plot 109 71 109 93 /

\plot 68 10 132 10 /
\plot 68 20 132 20 /
\plot 68 30 132 30 /
\plot 68 40 132 40 /
\plot 68 50 132 50 /

\plot 61 57 125 57 /
\plot 54 64 118 64 /
\plot 47 71 111 71 /

\plot 71.4 48.6 47.6 72.4 /
\plot 81.4 48.6 57.6 72.4 /
\plot 91.4 48.6 67.6 72.4 /
\plot 101.4 48.6 77.6 72.4 /
\plot 111.4 48.6 87.6 72.4 /
\plot 121.4 48.6 97.6 72.4 /
\plot 131.4 48.6 107.6 72.4 /

\plot 71.4 38.6 47.6 62.4 /
\plot 71.4 28.6 47.6 52.4 /
\plot 71.4 18.6 47.6 42.4 /
\plot 71.4 8.6 47.6 32.4 /
\put {\it Figure~5} [l] at 80 0
\endpicture

\vskip -25pt\hskip100pt
\beginpicture
\setcoordinatesystem units <.3mm,.3mm> point at 0 0
\setlinear
\plot 20 20 150 20 /
\put{$\bullet$} at 20 20

\plot 30 20 30 25 /
\plot 26.5 23.5 26.5 28.5 /
\plot 33.5 23.5 33.5 28.5 /
\plot 30 20 26.5 23.5 /
\plot 30 20 33.5 23.5 /
\plot 30 25 26.5 28.5 /
\plot 30 25 33.5 28.5 /
\plot 26.5 28.5 30 32 /
\plot 33.5 28.5 30 32 /

\plot 60 20 60 30 /
\plot 53 27 53 37 /
\plot 67 27 67 37 /
\plot 60 20 53 27 /
\plot 60 25 53 32 /
\plot 60 30 53 37 /
\plot 63.5 33.5 56.5 40.5 /
\plot 67 37 60 44 /
\plot 60 20 67 27 /
\plot 60 25 67 32 /
\plot 60 30 67 37 /
\plot 56.5 33.5 63.5 40.5 /
\plot 53 37 60 44 /
\plot 56.5 23.5 56.5 33.5 /
\plot 63.5 23.5 63.5 33.5 /

\plot 90 20 90 35 /
\plot 86.5 23.5 86.5 38.5 /
\plot 83 27 83 42 /
\plot 79.5 30.5 79.5 45.5 /
\plot 93.5 23.5 93.5 38.5 /
\plot 97 27 97 42 /
\plot 100.5 30.5 100.5 45.5 /
\plot 90 20 79.5 30.5 /
\plot 90 25 79.5 35.5 /
\plot 90 30 79.5 40.5 /
\plot 90 35 79.5 45.5 /
\plot 93.5 38.5 83 50 /
\plot 97 42 86.5 53.5 /
\plot 100.5 45.5 90 57 /
\plot 90 20 100.5 30.5 /
\plot 90 25 100.5 35.5 /
\plot 90 30 100.5 40.5 /
\plot 90 35 100.5 45.5 /
\plot 86.5 38.5 96.5 50 /
\plot 83 42 93.5 53.5 /
\plot 79.5 45.5 90 57 /

\plot 130 20 130 40 /
\plot 126.5 23.5 126.5 43.5 /
\plot 123 27 123 47 /
\plot 119.5 30.5 119.5 50.5 /
\plot 116 34 116 54 /
\plot 133.5 23.5 133.5 43.5 /
\plot 137 27 137 47 /
\plot 140.5 30.5 140.5 50.5 /
\plot 144 34 144 54 /
\plot 130 20 116 34 /
\plot 130 25 116 39 /
\plot 130 30 116 44 /
\plot 130 35 116 49 /
\plot 130 40 116 54 /
\plot 133.5 43.5 119.5 57.5 /
\plot 137 47 123 61 /
\plot 140.5 50.5 126.5 64.5 /
\plot 144 54 130 68 /
\plot 130 20 144 34 /
\plot 130 25 144 39 /
\plot 130 30 144 44 /
\plot 130 35 144 49 /
\plot 130 40 144 54 /
\plot 126.5 43.5 140.5 57.5 /
\plot 123 47 137 61 /
\plot 119.5 50.5 133.5 64.5 /
\plot 116 54 130 68 /
\put {\it Figure~6} [l] at 60 0
\setdashes
\plot 150 20  170 20 /

\endpicture

 \bigskip

The following example was suggested by D.~Cassi, R.~Burioni and A.~Vezzani.
\example{cubes}
{Let $X$ be the graph obtained by attaching at each vertex $i$ of $\naturali$
a cube lattice $C_i$ of side $n_i$ by one of its corners (in Figure 7 $n_i=i$).
Suppose that
$n_i$ diverges. Then the simple random walk on $X$ is locally
recurrent (by an application of
Theorem~\lemmaref{flow}), hence it is also \roa\ with respect
to any $\lambda$.
\hfill\break\indent
Consider the following ICF: $\cF=\{\bigcup_{i=1}^nC_i\}_n$. The simple random
walk on $X$ is \ttoa\ with respect to the limit on the average $L_\cF$.
Indeed for each $k\in\naturali$, the set $X_k$ obtained removing from $X$ all
the vertices at distance $k$ from the surface of the cubes has $L_\cF$-measure
equal to $1$ and if $x\in X_k$, $f^{(k)}(x,x)=f_{\zzz^3}^{(k)}$.
The thesis is an easy consequence of  Theorem~\lemmaref{Alexandroff1}
and Theorem~\lemmaref{limit1}(a.1).
Note that since the classification on the average and the thermodynamical one
are different this also provides an example of a random
walk for which $F(x,x)$ is not totally convergent
(Proposition~\lemmaref{twodef}(iv)).  }

\autosez{convexrw}{The general case}

In this section we consider a general sequence $\lambda=\{\lambda_n\}_n$
of probability measures on $X$. We already stated many results
for the average over balls which hold also for a general $\lambda$:
Proposition~\lemmaref{twodefo},
 Proposition~\lemmaref{Ginfinito}, Theorem~\lemmaref{flusso-o},
Proposition~\lemmaref{inducedo},
Corollary~\lemmaref{subtoa}, Theorem~\lemmaref{subtrano} and
Corollary~\lemmaref{subtran1o}.

We now prove a result which we already used in the previous sections:
it is a comparison between classifications with two different
limits on the average (see Propositions~\lemmaref{indip} and \lemmaref{FBC1}).

\proposition{twolimits}
{Let $\lambda=\{\lambda_n\}_{n}$,
$\eta=\{\eta_n\}_{n}$ two
sequences of probability measures on $X$.
Let us   consider the following assertions:
\vskip 4 pt
\item{(i)} there exist two divergent sequences $\{i_n\}_n$ and $\{j_n\}_n$
of natural numbers, and two positive constants  $C,K$ such that
$C\lambda_{i_n}(x)\ge\eta_{n}(x)$ and $K\eta_{j_n}(x)\ge\lambda_{n}(x)$
for every $n$ and $x$;
\item{(ii)} for every $A \subseteq X$, $A \in {\cal D}(L_\lambda)$,
$L_\lambda(A)=1$ if and only if
$A \in {\cal D}(L_\eta)$,
$L_\eta(A)=1$;
\item{(iii)} a random walk $(X,P)$
is \roa\  (respectively \toa) with respect to $\lambda$ if and
only if it is \roa\  (respectively \toa) with respect to
$\eta$. \vskip 4 pt
Then the following chain of implications holds:
(i) $\Ra$ (ii) $\Ra$ (iii).
}
\proof
(i) $\Ra$ (ii) It is easy (recall also Remark~\lemmaref{infsup1}).

(ii) $\Ra$ (iii) It is true because of the equivalence
between (ii) and (iii) in Proposition~\lemmaref{Alexandroff4}.
\QED

We state and prove for the general case some of the results
quoted in Sections~\sref{average} and \sref{sec:subgraphs}.

\proposition{twodef}
{Let $(X,P)$ be a random walk  and
let $\infty$ be the point added to $X$ in order to construct
its one point compactification.
\item{(i)} If there exists $A \subseteq X$ measurable, such
that $L_\lambda(A)=1$ and $\lim_{\scriptstyle x \rightarrow \infty\atop
\scriptstyle x\in A} F(x,x) =\alpha$ then
$L_\lambda(F)$  exists, and is equal to $\alpha$. Thus
the random walk is $\lambda$-\toa\ (respectively $\lambda$-\roa)
if and only if $\alpha < 1$ (respectively $\alpha=1$);
\item{(ii)} if $(X,P)$ is (locally) recurrent then
 $L_\lambda(F)$ exists, is equal to $1$ and the random walk is
$\lambda$-\roa;
\item{(iii)} if $(X,P)$ is $\lambda$-\troa\ then $L_\lambda(F)$ exists,
 is equal to $1$ and the random walk is $\lambda$-\roa;
\item{(iv)} if the series $F(x,x)$ is totally convergent (with
respect to $x\in X$) and $(X,P)$ is $\lambda$-\ttoa\ then $L_\lambda(F)$ exists,
it is less than $1$ and the random walk is $\lambda$-\toa;
\item{(v)} $(X,P)$ is $\lambda$-\roa\ $\Longleftrightarrow$ for every
 $\eps>0$ the set $\{x:F(x,x)\ge 1-\eps\}$ is measurable with measure $1$;
\item{(vi)}  if there exists
$A \subseteq X$ measurable, such that $L_\lambda(A)=1$ and
$\lim_{\scriptstyle x \rightarrow \infty\atop
\scriptstyle x\in A} F(x,x)=1$ then
$(X,P)$ is $\lambda$-\roa. Moreover, if $\lambda$ is regular the
converse holds;
\item{(vii)} $(X,P)$ is $\lambda$-\toa\ $\Longleftrightarrow$ there exists $A\sbs X$
  such that $\supL_\lambda(A)>0$ and $\sup_A F(x,x)<1$.
\item{}
\vskip -10 pt}

\proof

(i) It is an easy consequence of Proposition~\lemmaref{split1}.

(ii) It follows trivially from (i).

(iii) It is a consequence of Remark~\lemmaref{oss1} since $z \mapsto
F(x,x|z)$ is a non decreasing function on $[0,1]$ bounded from above by $1$.

(iv) It follows by Theorem~\lemmaref{limit1}(c).

(v) and (vii) From the relation $0 \leq \infL_\lambda(F) \leq 1$ we
have
 that $(X,P)$ is $\lambda$-\toa\ if and only if it is not $\lambda$-\roa;
Proposition~\lemmaref{Alexandroff4} yields the conclusion.

(vi) It is a consequence of Proposition~\lemmaref{Alexandroff4}.
\QED

\theorem{flusso}
{Let $(X,P)$ be a reversible random walk, with
reversibility
measure $m$ satisfying $\inf m(x)>0$, $\sup m(x)<+\infty$ (in particular
this condition is satisfied by the simple random walk on a graph with
bounded geometry).
Then TFAE:
\item{(a)} the random walk is $\lambda$-\toa;
\item{(b)} there exists $A\sbs X$ such that $\sup L_\lambda(A)>0$, there
is a finite  energy flow $u^x$ from $x$ to $\infty$ with non-zero input
$i_0$ for
every $x\in A$  and $\sup_{x\in A} <\!u^x,u^x\!><+\infty$;
\item{(c)} there exists $A\sbs X$ such that $\sup L_\lambda(A)>0$
 and $\inf_{x\in A} \Cap (x)>0$.}

\proof
First note that the interesting case is $(X,P)$ (locally) transient.
We therefore restrict to this particular case.

(a) $\Ra$ (b) Recall that $u^x=-{i_0\over m(x)}\nabla G(\cdot,x)$ is
a finite energy flow from $x$ to $\infty$ with input $i_0$ and energy
$$
<u^x,u^x>={i_0^2\over m(x)}\,G(x,x),
\autoeqno{ux}
$$
(where $\nabla$ denotes the difference operator, see \cite{Woess2}).
But by Proposition~\lemmaref{twodef}(vii) the network is $\lambda$-\toa\ if and only
if there exists $\alpha<1$,
$A\sbs X$ such that $\inf L_\lambda (A)>0$ and
$F(x,x)<\alpha<1$ for every $x\in A$.
Since $G(x,x)=1/(1-F(x,x))$ this is equivalent to
$\sup_{x\in A} G(x,x)<+\infty$. By equation \eqref{ux} and our hypotheses
on the reversibility measure, this implies (b).

(b) $\Ra$ (c) This is an obvious consequence of
$\Cap (x)\ge 1/<u^x,u^x>$ (see for instance 
\cite{Woess2}).

(c) $\Ra$ (a) This follows from $G(x,x)\le m(x)/\Cap(x)$ for every
$x\in A$, and from our hypotheses on $m$.
\QED

\theorem{subtran}
{Let $(A,E(A))$ be a subgraph of $X$ such that
$\supL_\lambda(A)>0$.  Suppose that
there exists $x_0\in A$ such that for every vertex $y \in A$ there
exists an injective map ${\gamma_y} : A \rightarrow A$ such that
(i) ${\gamma_y}(x_0)=y$ and  (ii) for any $w,z \in A$, $(w,z) \in E(A)$
implies
$({\gamma_y}(w), {\gamma_y}(z)) \in E(A)$. If the simple random walk on
$(A,E(A))$
is transient then the simple random walk on $(X,E(X))$ is $\lambda$-\toa.
}

\proof
Let us consider $(A,E(A))$ with the edge orientation induced by $X$.
Let $u$ be a flow on $(A,E(A))$ with finite energy starting from $x_0$ to
infinity with input $1$. Given any simple random walk, the conductance
is the characteristic function of the edges, hence for any $y\in A$
it is easy to show that the following equation
$$
u_y(a,b):=\cases{
\epsilon_{\gamma_y}(a,b) u({\gamma_y}^{-1}(a),{\gamma_y}^{-1}(b)) & if
        $(a,b) \in E(A)$; \cr
0 & if $(a,b) \not\in E(A)$;\cr}
\qquad \forall (a,b) \in E(X)
$$
(where $\epsilon_{\gamma_y}(a,b)$ is equal to $+1$ or $-1$ according to
$({\gamma_y}(a),{\gamma_y}(b))^+=(x,y)^+$ or not)
define a finite energy flow $u_y$ on $(X,E(X))$ starting from $y$ to $\infty$
with input one.
Apply now Theorem~\lemmaref{flusso}.
\QED

\proposition{Ginfinit}
{Let $(X,P)$ be an irreducible random walk.
If $F(\cdot,\cdot|z) \in {\cal D}(L_\lambda)$
for every $z \in (\eps,1)$ for some $\eps\in (0,1)$ 
and $\lim_{z \rightarrow 1^-}
L_\lambda(F(z)) =1
$
(respectively $L_\lambda(F) =1$)
then
$\lim_{z \rightarrow 1^-} \infL_\lambda(G(z))=+\infty$
(respectively $L_\lambda(G)=+\infty$).
}
\proof
From equation $G(x,x|z)=1/(1-F(x,x|z))$ (see 
\cite{Woess2})
 we have that for all $x \in X$
and for all $z \in \reali$, $|z| <1$,
$
G(x,x|z)= \phi(F(x,x|z)),
$
where $\phi(t):= 1 / (1-t)$.
By Jensen's inequality
$$
\phi \left (\sum_{x \in X}  F(x,x|z)\lambda_{n}(x) \right ) \leq
\sum_{x \in X} G(x,x|z)\lambda_{n}(x).
$$
If we take the limit as $n$ goes to infinity
of both sides of the previous equation, taking
into account the continuity
of $\phi$,
$$
\eqalign{
\phi(L_\lambda(F(z)))= & \lim_{n \rightarrow +\infty}
\phi \left (\sum_{x \in X} F(x,x|z)  \lambda_{n}(x) \right ) \cr
\leq  &
\liminf_{n \rightarrow +\infty}
\sum_{x \in X} G(x,x|z) \lambda_{n}(x) =:\infL_\lambda(G(z));  \cr
}
$$
hence
$$
\lim_{z \rightarrow 1^-} \infL_\lambda(G(z)) \geq
\liminf_{z \rightarrow 1^-}
\phi(L_\lambda(F(z))) =
+ \infty.
$$
The case $L_\lambda(F) =1$ is completely analogous (please note
that it could happen that $\sum_{x \in X}  G(x,x) \lambda_{n}(x)= +\infty$
for some $n \in \naturali$).
\QED

\semiautosez{A} Appendix A: limits on the average of general functions  \par

In this appendix we consider a very general setting: if not otherwise
stated $X$ is a countable set, $\lambda$ is a sequence of probability
measures on $X$ and $f$ is a 
function defined on $X$.
We look for sufficient conditions on $f$ for the existence of
$L_\lambda(f)$ and we study what can be said about its value.

First we observe that if $\lambda$ is regular, values taken by $f$
on finite subsets (which have $L_\lambda$-measure zero) do not influence
neither the existence nor the value of $L_\lambda(f)$.
In some sense only values  $f(x)$ for $x$ tending to $\infty$ matter,
where $\infty$ is the point at infinity of the one point compactification of $X$.

\theorem{Alexandroff1}
{Let $\lambda$ be regular, then
for any real valued function $f$,
$$
\liminf_{x \rightarrow \infty} f(x)\leq \infL_\lambda(f) \leq
\supL_\lambda(f) \leq \limsup_{x \rightarrow \infty} f(x).
$$
In particular if there exists $\lim_{x \rightarrow \infty} f =\alpha$
then $f\in \dom(L_\lambda)$ and $L_\lambda(f)=\alpha$.
Moreover, if $f$ is bounded and $A\sbs X$ such that $L_\lambda (A)=1$, then
$$
\liminf_{\scriptstyle{x \rightarrow \infty\atop x\in A}} f(x)\leq
\infL_\lambda(f) \leq   \supL_\lambda(f) \leq
\limsup_{\scriptstyle{x \rightarrow \infty\atop x\in A}}
f(x).   $$
}
\proof
We deal only with the first inequality. Suppose that
 $\liminf_{x \rightarrow \infty} f(x) = \alpha \not
= -\infty$ (otherwise there is nothing to prove).
Since
$$
\liminf_{x \rightarrow \infty} f(x) \equiv \sup_{S \subset X:
|S|<+\infty}
\inf_{x \not \in S} f(x),
$$
for every $\epsilon >0$
there exists a finite subset $S$ such that for every $x\not \in S$, $f(x) > \alpha -
\epsilon$ and  hence
$$
\sum_{x \in X}  f(x)\lambda_{n}(x)
\geq
\sum_{x \in S} f(x) \lambda_{n}(x)+
(1-\lambda_{n}(S))(\alpha-\epsilon)
{\buildrel n \rightarrow +\infty \over
\longrightarrow} \alpha -\epsilon,
$$
whence $\alpha\le \infL_\lambda(f)$.
The rest of the proof is analogous.
\QED

We note that the previous result means that, for $\lambda$ regular,
if $f(x)$ converges (``almost surely'') when $x$ goes to infinity then
$f\in \dom(L_\lambda)$  and $L_\lambda(f)$ does not
depend on the choice of $\lambda$ (only the notion of
``almost surely'' does).

Moreover, only the topological (discrete) structure of $X$ is involved;
if $X$ has a graph structure then the topology
is obviously independent of the choice of the edge set
(the topology is always the discrete one).

If we restrict to bounded functions, we get a stronger result.

\proposition{Alexandroff4}
{Let 
$f:X \rightarrow \reali$ be such that $ N \leq f(x) \leq
M$,
for every
$x\in X$). Then TFAE:\vskip 4 pt
\item{(i)} $f \in {\cal D}(L_\lambda)$ and $L_\lambda(f)=M$ (respectively $L_\lambda(f)=N$);
\vskip 4 pt
\item{(ii)} $\infL_\lambda(f)=M$ (respectively $\supL_\lambda(f)=N$);
\vskip 4 pt
\item{(iii)} $\forall \epsilon >0$, $L_\lambda(\{x:f(x)>M-\epsilon\})=1$
(respectively  $L_\lambda(\{x:f(x)<N+\epsilon\})=1$);
\vskip 4 pt
\item{(iv)} $\forall \epsilon >0$,
$L_\lambda(\{x:f(x)\leq M-\epsilon\})=0$ (respectively
$L_\lambda(\{x:f(x)\geq N+\epsilon\})=0$).
\vskip 4 pt
Moreover if $\lambda$ is regular
then \vskip 4 pt
\item{(v)} there exists a measurable set $A$ with measure $1$ such that
$\lim_{\scriptstyle x \rightarrow +\infty \atop x \in A} f(x) = M$  
(respectively $\lim_{\scriptstyle x \rightarrow +\infty \atop x \in A} f(x) = 
N$)  \vskip 4 pt  
is equivalent to each of the previous ones.
}
\proof
We consider the case involving the superior limit $M$ (the other one
is completely analogous). Let us
define   $F_\epsilon^+:=\{x:f(x)>M-\epsilon\}$ and $F^-_\epsilon := X
\setminus F_\epsilon^+$.   \smallskip
$(i) \Leftrightarrow (ii)$ and $(iii) \Leftrightarrow (iv)$ are trivial.
\smallskip
$(i) \Rightarrow (iii)$. For every $n \in \naturali$, $\epsilon >0$,
$
\sum_{x \in X} f(x) \lambda_n(x) \leq M
\lambda_n(F_\epsilon^+) +
(M-\epsilon)
\lambda_n(F_\epsilon^-),
$
whence
$$
0 \leq \supL_\lambda({F_\epsilon^-}) \leq {1 \over \epsilon}
(M-\infL_\lambda (f) ) \leq {1 \over \epsilon}
(M-L_\lambda (f) )=0.
$$
\smallskip
$(iii) \Rightarrow (i)$. By Remark~\lemmaref{zeroset},
$f\chi_{F_\epsilon^-} \in {\cal D}(L_\lambda)$
and $L_\lambda(f \chi_{F_\epsilon^-})=0$ for every $\epsilon>0$. Whence
for every $\epsilon >0$
$$
M \geq \supL_\lambda(f) \geq \infL_\lambda(f) = \infL_\lambda(f
\chi_{F_\epsilon^+})  \geq (M-\epsilon)L_\lambda({F_\epsilon^+}) =
M-\epsilon,  $$
which easily implies $(i)$.
\smallskip
$(v) \Rightarrow (iii)$.
Let $n \in \naturali$ and $A_n:=\{x \in X: f(x) > M-1/n\}$.
Since $A\setminus A_n$ is finite and
$A_n^c=(A\setminus A_n)\cup (A^c\setminus A_n)$, then
$L_\lambda(A_n^c)=0$.
\smallskip
$(i) \Rightarrow (v)$.
Let $\{B_n\}_n$ be an increasing sequence of finite subsets of $X$ such
that $\cup_{n\in\naturali} B_n=X$ (that is $\{B_n\}_n$ is a basis for the set
of  neighbours of $\infty$).
Let, for any $n\in \naturali$, $A_n$ defined as in the previous point
($A_n$ is non-empty since $L_\lambda( f) =M$).
Let us construct recursively
two increasing sequences $\{m_i\}_i$ and $\{n_i\}_i$ with values in
$\naturali$   satisfying
$$
\eqalign{
&
\lambda_m(A_i) > 1 -1/i, \qquad \forall m \geq m_i;
\cr  &
\lambda_m(A_i\cap B_{n_i})> 1 -1/i, \qquad
\forall m :   m_i\leq m< m_{i+1}. \cr
}
$$
This is possible since $\lim_{m\rightarrow +\infty}
 \lambda_m(A_i)=1$ for any $i\in \naturali$ and
since (using Monotone Convergence Theorem)
$\lim_{n\rightarrow +\infty}
\lambda_m(A_i \cap B_{n})=
 \lambda_m(A_i)>1-1/i
$ and the set $\{m: m_i \leq m<m_{i+1}\}$ is finite.
We prove now that $A:=\cup_{i=1}^\infty (A_i \cap B_{n_i})$ satisfies
the two conditions in $(v)$.
\hfill\break\noindent
By regularity we have that $L_\lambda(A\setminus B_{n_i})=1$ hence
$A\setminus B_{n_i} \not = \emptyset$ and $\infty$ is an
accumulation point for $A$. Moreover if $x\in A\setminus B_{n_i}$ we have that
$x\in A_j$ for some $j>i$ and hence $f(x) > M - 1/j > M -1/i$;
this proves that
$\lim_{\scriptstyle x \rightarrow +\infty \atop x \in A} f(x) = M$.
\hfill\break\noindent
If $m$ satisfies $m_i \leq m <m_{i+1}$ then
$$
\lambda_m(A) \geq
\lambda_m(A_i \cap B_{n_i}) > 1-1/i
$$
whence
$\lim_{m\rightarrow +\infty}\lambda_m(A) =1$.
\QED

In the hypotheses of the previous theorem, if $\{\lambda_n\}_n$ is regular and
if for some $\epsilon >0$, the set $\{x:f(x)>M-\epsilon\}$ (respectively
 $\{x:f(x)<N+\epsilon\}$) is finite, then
$L_\lambda(\{x:f(x)=M\}) =1$
(respectively $L_\lambda(\{x:f(x)=N\}) =1$).

\medskip
We show now that ${\cal D}(L_\lambda)$ is closed in the uniform
convergence topology and that $L_\lambda$ is continuous.

\proposition{uniformly1}
{Let
$(\Gamma,\geq)$ be a
directed, partially ordered set.
If $\{f_\gamma\}_\Gamma \subseteq {\cal D}(L_\lambda)$ is a net with
the property that $\lim_\gamma f_\gamma =:f$ holds uniformly with
respect to $x \in X$, then $f \in {\cal D}(L_\lambda)$ and
$
L_\lambda(f) = \lim_\gamma L_\lambda(f_\gamma)
$.
}
\proof
It is trivial to note that
Theorem~7.11 of \cite{Rudin2} holds considering a net of functions
instead of a sequence.

Since $f_\gamma - f \rightarrow 0$ in $l^\infty(X)$ and $f_\gamma \in {\cal
D}(L_\lambda)$ for every $\gamma \in \Gamma$, then it is easy to show
that  $\sum_{x \in X}f(x) \lambda_{n}(x) $ exists for every $n \in \naturali$
and
$$
g_\gamma(n):=\sum_{x \in X} f_\gamma(x) \lambda_{n}(x) {\buildrel
\gamma \over \longrightarrow} \sum_{x \in X} f(x)\lambda_{n}(x)
 =:\phi(n)
$$
uniformly with respect to $n \in \naturali$ (since $\sum_{x \in X}^\infty
\lambda_{n}(x)=1$).

Since $\lim_{n \rightarrow +\infty} g_\gamma(n) =: L_\lambda(f_\gamma)$,
we have that
$\lim_\gamma L_\lambda(f_\gamma)$ and $ L_\lambda(f):=
\lim_{n \rightarrow +\infty} \phi(n)$ both exist and they are equal.
\QED

The previous result leads to an alternative proof
Theorem~\lemmaref{Alexandroff1}:
it is enough to observe that if $o\in X$ is fixed and $f$ satisfies
$\lim_{x \rightarrow \infty} f(x)=\alpha$ then $f_n(x):= (f-\alpha) \cdot
\chi_{B(o,n)}+\alpha$
is a sequence of functions uniformly convergent to $f$,
such that $f_n \in{\cal D}(L_\lambda)$ and $L_\lambda(f_n)=\alpha$.


\remark{oss1}
{As in the case of the thermodynamical classification, $f$ may
depend not only on $x\in X$ but also on $z\in [0,r)$.
If $z \mapsto f(x,z)$ is non negative 
and  not decreasing as $z \in (\epsilon,r)$ (for each fixed $x$),
 then
$$
\lim_{z \rightarrow r^-}
\infL_\lambda(f(\cdot,z))=
\sup_{z \in (\epsilon,r)} \infL_\lambda (f(\cdot,z))\leq
\infL_\lambda(f(\cdot,r)),
$$
\smallskip
moreover if $\lim_{z \rightarrow r^-} \infL_\lambda (f(\cdot,z)) = \sup_{x
\in X} f(x,r) <+\infty$
then
$f(\cdot,r) \in{\cal D}(L_\lambda)$ and $L_\lambda(f(\cdot,r))=
\sup_{x \in X} f(x,r) <+\infty$.
}

\semiautosez{B} Appendix B: limits on the average of power series  \par

In this appendix $X$ and $\lambda$ are as in the
preceding section. We study the limit on the average of
families of power series.
In particular we search for conditions on
$\sum_{n=0}^\infty a_n(x) z^n$, to belong to $\dom(L_\lambda)$ for every
fixed $z$ in the common domain of convergence; moreover, provided that
$L_\lambda(\sum_{n=0}^\infty a_n(x) z^n)$ exists we ask when
$$
\lim_{z \rightarrow r^-}L_\lambda \left(\sum_{n=0}^\infty a_n(x) z^n \right)
=L_\lambda \left ( \sum_{n=0}^\infty a_n(x) r^n \right ),
$$
where $z \in \reali$ and $B(0,r)$ is a common domain of convergence.

\theorem{limit1}
{Let 
$\sum_{n=0}^\infty a_n(x)z^n$ be
a family of power series such that $a_n(x)\ge 0$ for every
$n\in \naturali, x \in X$. Suppose that the series
$\sum_{n=0}^\infty k_n z^n$, where $k_n:= \sup_{x \in X} a_n(x)$,
 has a positive radius of convergence $r^\prime$.
If $r \in (0, r^\prime]$ then the
following results hold: \smallskip
\item{(a)} if $\{a_n(\cdot)\}_n \subset \dom(L_\lambda)$ then
\itemitem{(a.1)} $\sum_{n=0}^\infty a_n(\cdot) z^n \in \dom(L_\lambda)$ and
$L_\lambda(\sum_{n=0}^\infty a_n(\cdot) z^n)=
\sum_{n=0}^\infty L_\lambda(a_n) z^n$,
for every $z\in \complessi$ such that
$\sum_{n\in \naturali} k_n |z|^n<+\infty$;
\smallskip
\itemitem{(a.2)} $\lim_{z \rightarrow r^-} L_\lambda(\sum_{n=0}^\infty
a_n(\cdot) z^n)= \sum_{n=0}^\infty L_\lambda(a_n) r^n
\leq \infL_\lambda(\sum_{n=0}^\infty a_n(\cdot) r^n)$
(the values possibly being infinite),
\smallskip
\itemitem{(a.3)} if $\sum_{n=0}^\infty k_n r^n<+\infty$ then
$\lim_{z \rightarrow r^-} L_\lambda(\sum_{n=0}^\infty
a_n(\cdot) z^n)=  L_\lambda(\sum_{n=0}^\infty a_n(\cdot) r^n)$; \smallskip
\item{(b)} if $\sum_{n=0}^\infty a_n(\cdot) z^n \in \dom(L_\lambda)$ for all
$z\in\Gamma$ where $\Gamma$ is a $C^1$-circuit
with $\sup_{z\in\Gamma}\vert z\vert =r_1<r$, and $d(\Gamma,0)=r_2>0$  then
$\{a_n\}\subset \dom (L_\lambda)$ and (a.1) holds;\smallskip
\item{(c)} if $\sum_{n=0}^\infty k_n r^n<+\infty$, then
   $\lim_{k\ra\infty} L_\lambda(\sum_{n=0}^\infty a_n(x) z_k^n)
     =L_\lambda(\sum_{n=0}^\infty a_n(x) r^n)$
  provided that $\vert z_k\vert<r$, $z_k\ra r$ and
  $\sum_{n=0}^\infty a_n(\cdot) z_k^n \in{\cal D}(L_\lambda)$ for
     every $k$.
\item{}
}

\proof
The proof is based on quite classical arguments,hence we just outline it.

{\it (a.1)} Exchange the order of summation (using for instance
Fubini-Tonelli's Theorem) and apply Bounded Convergence Theorem
to the series.

{\it (a.2)} It is an easy application of Monotone Convergence Theorem.

{\it (a.3)} It follows from (a.1) and (a.2).

{\it (b)} Using the Cauchy integral formula and Fubini's Theorem
we obtain
$$
\sum_{x\in X} a_k(x)\lambda_n(x) ={1\over 2\pi i} \int_\Gamma {\sum_{x \in X}
\left ( \sum_{j =0}^\infty a_j(x) z^j \lambda_n(x) \right )
\over z^{k+1}} \d z,
$$
and simple calculations show that the norm of
the integrand is bounded by a constant
(uniformly with respect to  $n\in \naturali$),
hence Bounded Convergence Theorem implies the result.

{\it (c)} It follows by    Proposition~\lemmaref{uniformly1},
since   $\lim_{\scriptstyle z\ra r\atop\scriptstyle \vert z\vert<r}
\sum_{n=0}^\infty a_n(x)z^n=
\sum_{n=0}^\infty a_n(x)r^n$ holds uniformly with respect to
$x\in X$.
\QED

The meaning of the previous theorem is that for a
family of power series satisfying the hypotheses
of the theorem,
the  following assertions are equivalent: \smallskip
\item{(i)} every coefficient is in the domain of $L_\lambda$; \smallskip
\item{(ii)} there exists a circuit $\Gamma$ as in Theorem~\lemmaref{limit1}(b),
such that for all $z\in\Gamma$,
$\sum_{n=0}^{+\infty} a_n(x) z^n$ is in the domain of $L_\lambda$;\smallskip
\item{(iii)} for every $z\in \complessi$, $|z| < r^\prime$,
$\sum_{n=0}^{+\infty} a_n(x) z^n$ is in the domain of $L_\lambda$.\smallskip

By means of the previous theorem we can state and prove a result
which we call {\it identity principle on the average for power series}.

\proposition{identity}
{If $w_j(x,z):=\sum_{i=0}^\infty a^{(j)}_i(x) z^i$,
$j=1,2$,
is a couple of families of power series on 
$X$, with non negative coefficients. Suppose that the series
$\sum_{i=0}^\infty k^j_i z^i$, where $k^j_i:= \sup_{x \in X} a_i^{(j)}(x)$,
 has a positive radius of convergence $r_j$, $j=1,2$,
and that $w_j(\cdot,z) \in {\cal D}(L_\lambda)$  $j=1,2$, for all $z \in
B(0,\min(r_1,r_2))$. Then TFAE \smallskip
\item{(i)} there exists a subset $E\subseteq B(0, \min(r_1,r_2))$
with an accumulation  point $x_0$ which belongs
to the domain $B(0,\min(r_1,r_2))$ such that   
$$ 
L_\lambda(w_1(\cdot, z))=L_\lambda(w_2(\cdot, z)), \quad \forall z \in E; 
$$ 
\item{(ii)} for every $n \in \naturali$ and for every $x\in X$ we have that 
$L_\lambda(a^{(1)}_n)= L_\lambda(a^{(2)}_n)$. 
} 
\proof 
By Theorem~\lemmaref{limit1}(a.1) we have that
$$ 
L_\lambda(w_j(\cdot,z))=\sum_{i=0}^\infty L_\lambda(a^{(j)}_i) z^i,\quad 
\forall z\in B(o,r_j), \ j=1,2, 
\autoeqno{a.1bis} 
$$
then (i) $\Longrightarrow$ (ii) is trivial.

\noindent (ii) $\Longrightarrow$ (i) It is a consequence of \eqref{a.1bis}
and Theorems~8.1.2 and 8.1.3 of \cite{Hille}.
\QED

We can state a similar result which takes also into account the presence
of zero-measure sets.

\proposition{identity2}
{Suppose that 
$w_i(x,z):=\sum_{j=0}^\infty a^{(i)}_j(x) z^j$, $i=1,2$
satisfy the hypotheses of the previous theorem, and that
for every $n\in \naturali$ there exists a subset $X_n \subset X$ such that
$L_\lambda({X_n}) =0$ and $a^{(1)}_n(x)=a^{(2)}_n(x)$ for every
$x \not \in X_n$.
If $z\in \complessi$, $|z| \in (0, \min(r_1,r_2))$ then
$$
w_1(\cdot, z)\in {\cal D}(L_\lambda)
\quad \Longleftrightarrow \quad
w_2(\cdot, z)\in {\cal D}(L_\lambda),
\autoeqno{equiv1}
$$
and
$$
L_\lambda(w_1(\cdot,z)) =L_\lambda(w_2(\cdot, z)).
\autoeqno{equiv2}
$$
}
\proof
Given any general family of power series
$w(x,z):=\sum_{i=0}^\infty a_i(x)z^i$ such that
$\sum_{i=0}^\infty k_iz^i$ has
a positive radius of convergence $r$ ($k_i:=\sup_{x \in X}
|a_i(x)|$), using Bounded Convergence Theorem
it is not difficult to prove that (uniformly with respect to $|z|\leq r-\epsilon$),
$$
\lim_{n \rightarrow +\infty} \sum_{i=o}^\infty \sum_{x\in X_i}
a_i(x)  z^i \lambda_{n}(x) = 0.
$$

Hence if $|z|<r$ it is obvious that
$$
\eqalign{
w_j(\cdot, z) \in {\cal D}(L_\lambda) \ \Longleftrightarrow \  \exists
\lim_{n \rightarrow +\infty} \sum_{i=0}^\infty \sum_{x \in X_i^c}
 a_i^{(j)}(x) z^i \lambda_{n}(x),\cr
L_\lambda(w_j(\cdot,z)) =
\lim_{n \rightarrow +\infty} \sum_{i=0}^\infty \sum_{x \in X_i^c}
a_i^{(j)}(x) z^i\lambda_{n}(x) , \cr
}
$$
but by our hypotheses $a_i^{(1)}(x)= a_i^{(2)}(x)$ for every $i \in \naturali$ and
for every $x\in X_i^c$, whence the existence and the value of the last limit does not
depend on $j$.
\QED


\autosez{conclusion}{Conclusions}

In this paper we proposed a new classification
and we showed how is it possible to manage it from a 
technical point of view. However, following \cite{Cassi},
it seems reasonable to consider also the ``dual'' classification
obtained by substitution in definition \lemmaref{convexlimit} 
of the $\liminf$ with the $\limsup$.
With slight differences, it is possible to prove similar
results for this new classification. 
Here we state two crucial results (Proposition~
\lemmaref{Alexandroff5} and Theorem~\lemmaref{supflusso})
and we give the proof of the former, the proof of the latter
being scarcely different from its analogous.

\definition{Suprec}{Let $(X,P)$ be a random walk
and $\{\lambda_n\}$ a sequence of probability measures
on $X$. The random walk is called {\it $\lambda$-Suprecurrent}
(resp.~{\it $\lambda$-Suptransient}) if and only if
$\supL_\lambda(F)=1$
(resp.~$\supL_\lambda(F)<1$). 
}

We immediately note that $\lambda$-\roa implies
$\lambda$-Suprecurrent; moreover a random walk
is $\lambda$-Suptransient if and only if 
there exist $\epsilon>0$ and $n_0 \in \naturali$ such that
for any $n \geq n_0$ we have $\sum_{x\in X} F(x) \lambda_n(x) <1-\epsilon$.

The following proposition the analogous of
Proposition~\lemmaref{Alexandroff4}.

\proposition{Alexandroff5}{Let $F:X \rightarrow \reali$ such that
$N \leq F(x) \leq M$ for any $x\in X$; if $\lambda$ is regular then
consider the following assertions
\item{(i)} $\supL_\lambda(F)=M$ (resp.~$\infL_\lambda(F)=N$);
\item{(ii)} there exists a subset $A \subseteq X$ such that $\supL_\lambda(A)=1$ and
$\lim_{\scriptstyle x \rightarrow +\infty \atop \scriptstyle x \in A}
F(x) =M$
(resp.~$\lim_{\scriptstyle x \rightarrow +\infty \atop \scriptstyle x \in A}
F(x) =N$);
\item{(iii)} for any  $\epsilon >0$, $\supL_\lambda(\{F>M-\epsilon\})=1$
($\supL_\lambda(\{F>N+\epsilon\})=1$).
Hence $(i) \Longleftrightarrow (iii)$; moreover if $\lambda$ is regular
$(ii)$ is equivalent to any of the previous ones.}
\proof
$(i) \Longrightarrow (ii)$. Let $\{n_j\}_{j \in \naturali}$ be such that
$\lim_{j \rightarrow +\infty} \sum_{x\in X}F(x) \lambda_{n_j}(x) =M$
and let $\eta_j:=\lambda_{n_j}$ for any $j \in \naturali$.
Then $L_\eta(F)=M$, hence, by Proposition~\lemmaref{Alexandroff4}
(being $\eta$ regular), there exists $A\subseteq X$ such that
$L_\eta(A)=1$ and 
$\lim_{\scriptstyle x \rightarrow +\infty \atop \scriptstyle x \in A}
F(x) =M$ (obviously
$L_\eta(F)=1$ implies
$\supL_\lambda(F)=1$).

$(ii) \Longrightarrow (i)$. Let $\{n_j\}_{j \in \naturali}$
be such that $\lim_{j \rightarrow +\infty} \lambda_{n_j}(A)=1$
and, given any $\epsilon >0$, let $K_\epsilon \subseteq A$ 
such that $|A \setminus K_\epsilon|<+\infty$ and 
$F|_{K_\epsilon} >M-\epsilon$. It is easy to show that
$$
\eqalign{
\sum_{x \in X} F(x) \lambda_n(x) &= \sum_{x\in A} F(x) \lambda_n(x)+
\sum_{x\in A^c} F(x) \lambda_n(x) \geq{} \cr
{} &\geq \lambda_n(K_\epsilon)(M-\epsilon)+
\lambda_n(A^c) N;
}
$$
hence 
$$
M \geq \limsup_{j \rightarrow +\infty} \sum_{x \in X} F(x) \lambda_{n_j}(x)
\geq \limsup_{j \rightarrow +\infty} \lambda_{n_j}
(K_\epsilon)(M-\epsilon)+
\limsup_{j \rightarrow +\infty} \lambda_{n_j}(A^c) N= M-\epsilon
$$
since 
$\lambda_n(K_\epsilon)=\lambda_n(A)-\lambda_n(A\setminus K_\epsilon)$
and, by regularity, $\lim_{j \rightarrow +\infty} \lambda_{n_j} 
(A \setminus K_\epsilon) =0$.

$(i) \Longrightarrow (iii)$. Let $\eta_k:=\lambda_{n_k}$ where
$\{n_k\}_k$ is such that 
$$
\lim_{k\rightarrow +\infty} 
\sum_{x \in X} F(x) \lambda_{n_k} (x) =M,
$$
then, for every $\epsilon >0$,
$$
1 \geq \supL_\lambda(\{F>M-\epsilon\}) 
\geq L_\eta(\{F>M-\epsilon\}) =1
$$
since we apply Proposition~\lemmaref{Alexandroff4}
to $\eta$.

$(iii) \Longrightarrow (i)$. We know that there exists
a non decreasing sequence $\{n_k\}_k$ such that 
for any $k>0$
$$
\sum_{x \in X} \chi_{\{F>M-1/k\}}(x) \lambda_n(x)>1 - {1 \over k},
\quad \forall n \geq n_k.
$$
If $\eta_k:=\lambda_{n_k}$ then, given any
$\epsilon >0$ and $\delta >0$,
for $k \geq 1/\epsilon+1/\delta$ we have that
$$
\sum_{x \in X} \chi_{\{F>M-\epsilon\}}(x) \lambda_n(x)
\geq
\sum_{x \in X} \chi_{\{F>M-1/k\}}(x) \lambda_n(x)
\geq 1-{1\over k} > 1 -\delta,
$$
hence
$L_\eta(\{F>1-1/\epsilon\}) =1$ and,
by Proposition~\lemmaref{Alexandroff4},
there exists $A \subseteq X$ such that
$L_\eta(A)=1$ and 
$\lim_{\scriptstyle x \rightarrow +\infty \atop \scriptstyle x \in A}
F(x) =M$ which implies
$\supL_\lambda(A)=1$
\QED

Condition (ii) of the previous theorem is
technically easy and it allows us, in the 
natural case, to characterize suprecurrence.
The following theorem is the analogous of
Proposition~\lemmaref{flusso}: the proof
is omitted since it is now straightforward.

\theorem{supflusso}
{Let $(X,P)$ be a reversible random walk, with
reversibility
measure $m$ satisfying $\inf m(x)>0$, $\sup m(x)<+\infty$ (in particular
this condition is satisfied by the simple random walk on a graph with
bounded geometry).
Then TFAE:
\item{(a)} the random walk is $\lambda$-Suptransient;
\item{(b)} there exists $A\sbs X$ such that $\infL_\lambda(A)>0$, there
is a finite  energy flow $u^x$ from $x$ to $\infty$ with non-zero input
$i_0$ for
every $x\in A$  and $\sup_{x\in A} <\!u^x,u^x\!><+\infty$;
\item{(c)} there exists $A\sbs X$ such that $\infL_\lambda(A)>0$
 and $\inf_{x\in A} \Cap (x)>0$.}

\vskip18pt\noindent
{\bf Acknowledgments}
  \smallskip

We would like to thank R.~Burioni, D.~Cassi and A.~Vezzani for the stimulating
conversations.

   The authors acknowledge support from the E.S.I. (Wien, Austria).
The first author acknowledges also support from the ARGE Alpen-Adria
(Steiermark - Austria).

\vskip 40 pt  
{\bf Bibliography}  
\bigskip  
\insertbibliografia  
  
\end